\documentclass[hidelinks,onefignum,onetabnum]{siamart220329}

\usepackage{amsfonts,amsmath}
\usepackage{graphicx}
\usepackage{epstopdf}
\usepackage{algorithmic}
\usepackage{subfig,multicol,dashrule}

\definecolor{Green}{RGB}{53, 181, 53}

\newsiamremark{remark}{Remark}
\newsiamremark{assumption}{Assumption}

\newcommand*{\hl}{}

\newcommand*{\tm}{}
\def\hlb{\color{black}}
\headers{Infinite-Horizon Discrete-Time Dynamic Games}{A. Monti, B. Nortmann, T. Mylvaganam, M. Sassano}

\title{Feedback and Open-Loop Nash Equilibria for LQ Infinite-Horizon Discrete-Time Dynamic Games\thanks{
\funding{The work of B. Nortmann is supported by UKRI DTP 2019 grant number EP/R513052/1.
		The work of M. Sassano is partially supported by the Italian Ministry for Research in the framework of the 2020 Program for Research Projects of National Interest (PRIN). Grant No. 2020RTWES4.}}}

\author{Andrea Monti\thanks{Institute for Software Technology at the German Aerospace Center (DLR), 38108 Braunschweig, Germany. This work has been carried out while he was a Ph.D. student at the Dipartimento di Ingegneria Civile ed Ingegneria Informatica (DICII), Università degli Studi di Roma ``Tor Vegata'', Via del Politecnico 1, 00133, Roma, Italy (\email{andrea.monti@uniroma2.it}).}
\and Benita Nortmann\thanks{Department of Aeronautics, Imperial College London, London SW7 2AZ, UK (\email{benita.nortmann15@imperial.ac.uk,t.mylvaganam@imperial.ac.uk}).}
\and Thulasi Mylvaganam\footnotemark[3]
\and Mario Sassano\thanks{Dipartimento di Ingengeria Civile e Ingegneria Informatica, Università di Roma, ``Tor Vergata''
(\email{mario.sassano@uniroma2.it}).}
}

\begin{document}

\maketitle

\begin{abstract}
We consider dynamic games defined over an \emph{infinite horizon}, characterized by linear, discrete-time dynamics and quadratic cost functionals. Considering such linear-quadratic (LQ) dynamic games, we  
focus on their solutions in terms \hl{of} Nash equilibrium strategies. Both Feedback (F-NE) and Open-Loop (OL-NE) Nash equilibrium solutions are considered. The contributions of the paper are threefold. First, our detailed study reveals some interesting structural insights in relation to F-NE solutions. Second, as a stepping stone towards our consideration of OL-NE strategies, we consider a specific infinite-horizon discrete-time (single-player) optimal control problem, wherein the dynamics are influenced by a known exogenous input and draw connections between its solution obtained via Dynamic Programming and Pontryagin's Minimum Principle. Finally, we exploit the latter result to provide a characterization of OL-NE strategies of the class of infinite-horizon dynamic games. The  results and key observations made throughout the paper are illustrated via a numerical example.
\end{abstract}

\begin{keywords}
Infinite-horizon Dynamic Games; Dynamic Programming; Pontryagin's Minimum Principle.
\end{keywords}


\section{Introduction} \label{sec:intro}
Game theory provides a framework to study strategic interactions between multiple \emph{decision makers}, referred to as \emph{players}, which may (or may not) be in competition with one another.
\emph{Dynamic} games, in particular, concern such \emph{multi-player decisions} that evolve over time via a dynamical system whose evolution is \emph{steered} by inputs of a certain number of \emph{players}. 
Such problems have gained interest across the control engineering field in the last decades and find a wide range of applications, for instance, in the context of cyber-security (\cite{huang2020,basar2019}), economics (\cite{jorgensen2010,long2011,jaskiewicz2022}), epidemiology (\cite{kordonis2022,nadini2020}), robotics (\cite{Myl2017,Cappello2021, Li-Bur:19, Nor-Myl:ECC22}) and power systems
	(\cite{nagata2002,pratap2017}).  Notably, dynamic games have gained increasing interest in recent years, in part due to the widespread presence of \emph{multi-agent systems} in modern applications, namely systems consisting of individual decision makers (or players) that interact, potentially in a competitive way, with one another (see \emph{e.g} \cite{chen2019control} and references therein). 
The framework provided by dynamic games has also proved useful for a wider range of problems involving dynamical systems with multiple inputs and/or multiple (potentially conflicting) objectives, such as robust control, mixed $H_2/H_\infty$ control and distributed control  (\cite{sweriduk1994differential,limebeer1992mixed,basar2008h,jaskiewicz2011stochastic,Cappello2022, Wrz20}).

In this paper we study \emph{infinite-horizon, dynamic} games.
More precisely, we consider a discrete-time linear system, the evolution of which is influenced by the actions (\emph{i.e.} control inputs) of a set of players. Each player aims at minimizing, via an individually manipulated input, its own quadratic cost functional over an infinite-horizon. 
While dynamic games can be considered a generalization of optimal control problems (which can be interpreted as a single-player game),
in the case of games the outcome for one player is dependent on the strategies of all other players. Hence, several alternative solution concepts exist for dynamic games. In what follows we focus on one of the most common notions of solution, namely \emph{Nash Equilibrium} (NE) strategies. These can be further categorized according to the information pattern available to each player while implementing its own strategy: \emph{Feedback Nash Equilibrium} (F-NE) 
and \emph{Open-Loop Nash Equilibrium} (OL-NE). 
See, \emph{e.g.} \cite{basar1998, starr1969further, Starr:69}, 
for further details on dynamic games and their various solution concepts.
Even in the case of linear-quadratic (LQ) games considered herein, NE solutions may not be unique, and the equations characterizing solutions are generally challenging to solve.
Somewhat surprisingly, the class of \emph{infinite}-horizon LQ \emph{discrete-time} dynamic games considered herein has received relatively limited attention in the literature hitherto (see, \emph{e.g.} \cite{basar1998} and \cite{freiling1999discrete}) compared to their continuous-time counterpart, namely \emph{differential} games defined over an infinite horizon. In the continuous-time setting infinite-horizon \emph{differential} games have been studied extensively (see, \emph{e.g.} \cite{Engwerda2005} and references therein). In this setting it is known that the application of Dynamic Programming (DP) arguments leads to a set of coupled Algebraic Riccati equations (AREs) characterizing F-NEs, whereas
arguments borrowed from Pontryagin's Minimum Principle (PMP) lead to (deceptively similar) \emph{asymmetric} AREs characterizing OL-NEs. 

The main contribution of this work is a complete characterization of all OL-NE solutions of \emph{infinite-horizon}, \emph{discrete-time} LQ dynamic games by relying on insights gained from PMP. Moreover, we provide conditions in the form of asymmetric coupled matrix equations (via invariance arguments) under which OL-NE strategies admit a \emph{feedback synthesis}.
As such, one of the main contributions of this paper can be seen as the discrete-time counterpart of the results presented in Chapters 7 and 8 of\footnote{Note that the terminology used throughout this paper is deliberately consistent with that of \cite{Engwerda2005} and the interested reader is referred to \cite{Engwerda2005} for further insights on LQ \emph{differential} games.} \cite{Engwerda2005}. 
Obtaining these results requires in turn 
the study of a special class of discrete-time optimal control problems, over an infinite horizon, wherein the system dynamics are influenced by a \emph{known} exogenous input. While this problem has been addressed in \cite{engwerda_infinityHor} via DP arguments, we draw novel connections with PMP by relating the solution to a certain trajectory of the associated Hamiltonian dynamical system. These novel insights are instrumental towards obtaining the aforementioned results related to OL-NE solutions. 
In addition, we provide a comprehensive characterization of F-NE in LQ dynamic games via DP arguments. Although similar conditions have appeared in the literature (see, \emph{e.g.} \cite[Sec. 6.2.3]{basar1998}), the specific matrix equations derived herein reveal an interesting structure, which allows comparison to both the equations characterizing OL-NE strategies and the continuous-time equivalent. Furthermore, we provide sufficient structural conditions (easily verifiable based solely on the problem data and more general than those in \cite[Sec. 6.2.3]{basar1998}) implying the asymptotic stability of the origin for the plant in closed loop with the F-NE strategies, provided they exist.
For ease of exposition, we present the analysis for two-player dynamic games.
However, the results are readily generalized to the $N$-player setting (at the cost of more cumbersome notation).

The remainder of the manuscript is organized as follows. In Section \ref{sec:problem} we recall some preliminaries on LQ dynamic games and their solutions in terms of NE strategies. F-NE strategies for infinite-horizon dynamic games are discussed in Section \ref{sec:F_NE}, while OL-NE strategies are derived 
in Section \ref{sec:OL_NE}. Finally, an illustrative numerical example is presented in Section~\ref{sec:sim}, before some concluding remarks are provided in Section~\ref{sec:conc}. 

\textbf{Notation} The set of real numbers \tm{and the set of non-negative natural numbers are} denoted by $\mathbb{R}$ and $\mathbb{N}_{\ge0}$, respectively. 
Given a square matrix $A$, $\sigma(A)$ indicates the set of the eigenvalues of $A$, namely the \emph{spectrum} of $A$, while $A^{-\top}$ represents the transpose of the inverse of $A$, namely $A^\top A^{-\top}~=~I$, and $A^\mathsf{H}$ denotes the conjugate transpose of $A$.  Given a symmetric matrix $A$, $A\succ 0$ ($A \succeq$ 0) denotes that the matrix is positive definite (positive semi-definite). Given a non-square matrix $B$ of full column rank, $B^\dagger$ denotes its 
(left) pseudo-inverse.  The observability matrix of the pair of matrices $(A,C)$ is denoted by $\mathcal{O}(A,C)$. $\mathbb{C}_{d}$ denotes the set $\{ \lambda \in \mathbb{C} : |\lambda| < 1 \}$.
The space of square summable infinite sequences $w = (w_0, w_1, w_2, ...)$, \emph{i.e.} such that $\sum_{k = 0}^{\infty} |w_k|^2 < \infty$, is denoted by $\ell_2$.

\section{Problem statement and preliminaries}
\label{sec:problem}
The main objective of this section consists in briefly recalling a few definitions and results concerning the computation of NE strategies in the setting of infinite-horizon LQ \hl{dynamic} games. 
Thus, consider a linear system described by the difference equation
\begin{equation}\label{eq:lti}
	x(k+1) = Ax(k) + B_1u_1(k) + B_2u_2(k),
\end{equation}
where $x : \tm{\mathbb{N}_{\ge0}} \rightarrow \mathbb{R}^{n}$ denotes the state of the system and $u_i : \tm{\mathbb{N}_{\ge0}} \rightarrow \mathbb{R}^{m_i}$, for $i = 1, 2$, denote the control actions of players 1 and 2, respectively, subject to the initial condition $x(0) = x_0 \in \mathbb{R}^{n}$. The aim of each player is to minimize an individual cost functional defined over an infinite horizon, described by
\begin{equation}
	\label{eq:cost_J_i}
	J_i(x_0, u_1(\cdot), u_2(\cdot)) = \frac{1}{2}\sum\limits_{k=0}^{\infty} x(k)^\top Q_ix(k) + u_i(k)^\top R_iu_i(k),
\end{equation}
with $Q_i = Q_i^{\top} \succeq 0$ and $R_i = R_i^{\top} \succ 0$ parameterized with respect to the initial condition $x_0$. \hl{In the game described by \eqref{eq:lti}, \eqref{eq:cost_J_i} (as well as in all, single or multi-player, dynamic optimization problems arising throughout the manuscript) it is further required that feasible solutions are such that the state tends to the origin as time diverges to infinity.} Since the value of the cost functional of each player depends
on the action selected by its opponent, indirectly via the evolution of the state,  
the cost functionals
can seldom be  minimized simultaneously. Thus the players need to settle for \emph{equilibrium}, rather than \emph{optimal}, solutions. 
Two common solution concepts that reflect different underlying information structures, namely feedback Nash equilibria (F-NE) and open-loop Nash equilibria (OL-NE),
	are recalled below (see, \emph{e.g.} \cite{basar1998}
	for further details).

\smallskip

\begin{definition}
	\label{def:FNE}
	(F-NE) A pair of control laws $(u_1^{\star}, u_2^{\star}) = (K_1^{\star}x, K_2^{\star}x)$ is a \emph{Feedback Nash Equilibrium} (F-NE) strategy if the inequalities
	\begin{subequations}\label{eq:def_FNE}
		\begin{align}
			J_1(x_0, K_1^{\star}x, K_2^{\star}x) & \leq J_1(x_0, K_1 x, K_2^{\star}x) \label{eq:def_FNE_ply_1}\,, \\[3pt]
			J_2(x_0, K_1^{\star}x, K_2^{\star}x) & \leq J_2(x_0, K_1^{\star} x, K_2x) \label{eq:def_FNE_ply_2}\,,
		\end{align}
	\end{subequations}
	hold for all $x_0 \in \mathbb{R}^{n}$ and for any pair $(K_1, K_2^{\star})$ such that $\sigma(A+B_1 K_1 + B_2 K_2^{\star}) \subset \mathbb{C}_d$ in the case of \eqref{eq:def_FNE_ply_1}, and for any pair $(K_1^{\star}, K_2)$ such that $\sigma(A+B_1 K_1^{\star} + B_2 K_2) \subset \mathbb{C}_d$ in the case of \eqref{eq:def_FNE_ply_2}.
\end{definition}

An appealing property of the solution concept detailed in Definition~\ref{def:FNE} is that the pair of strategies\footnote{As is common in the context of LQ dynamic games (see \emph{e.g.} \cite{Engwerda2007}), we limit our attention to linear feedback strategies in the context of F-NE solutions.} $(K_1^\star x, K_2^\star x)$ constitutes a F-NE solution of the game defined by the system~\eqref{eq:lti} and the cost functionals \eqref{eq:cost_J_i}, for $i = 1,2$, for \emph{all} initial conditions $x_0 \in \mathbb{R}^n$. Similarly, if at any point in time $\bar k > 0$ the actual state $x(\bar k)$ of the system~\eqref{eq:lti} deviates for any reason from the state $x^\star(\bar k)$, induced by the initial condition and the equilibrium strategies $(u_1^\star, u_2^\star)$ played for $k = 0, \ldots, \bar k - 1$, then the pair of strategies $(K_1^\star x, K_2^\star x)$ still constitutes a F-NE solution of the (restricted) game on $k \geq \bar k$. However, the feedback strategies entail 
that the players must measure the entire state of the underlying plant at all times in order to implement a F-NE strategy.
In practice, it may be interesting to study also the setting in which each player has knowledge of the initial condition $x_0$ only, together with the elapsed time, as discussed in the following definition. To streamline the presentation, let $\mathcal{U}_1(x_0) := \{ (\phi_1, \phi_2^{\star}) \in \ell_2 \times \ell_2 : J_1(x_0, \phi_1, \phi_2^{\star}) < \infty, \lim_{k \rightarrow \infty} x(k)=0 \}$ and $\mathcal{U}_2(x_0) := \{ (\phi_1^{\star}, \phi_2) \in \ell_2 \times \ell_2 : J_2(x_0, \phi_1^{\star}, \phi_2) < \infty, \lim_{k \rightarrow \infty} x(k) = 0 \}$.

\medskip 

\begin{definition}
	\label{def:OLNE}
	(OL-NE) Let $x_0 \in \mathbb{R}^{n}$ be given.
	A pair of control laws $(u_1^{\star}, u_2^{\star}) = (\phi_1^{\star}(k,x_0), \phi_2^{\star}(k,x_0))$ is an \emph{Open-Loop Nash Equilibrium} (OL-NE) strategy if the inequalities
	\begin{subequations}
		\label{eq:def_OLNE}
		\begin{align}
            J_1(x_0, \phi_1^{\star}, \phi_2^{\star})
            & \leq J_1(x_0, \phi_1, \phi_2^{\star})\,, \label{eq:def_OLNE_ply_1} \\[3pt]
            J_2(x_0, \phi_1^{\star}, \phi_2^{\star})
            & \leq J_2(x_0, \phi_1^{\star}, \phi_2)\,, 
            \label{eq:def_OLNE_ply_2}
		\end{align}
	\end{subequations}
	hold for any $(\phi_1, \phi_2^{\star}) \in \mathcal{U}_1(x_0)$, 
	in \eqref{eq:def_OLNE_ply_1} 
	and for any $(\phi_1^{\star}, \phi_2) \in \mathcal{U}_2(x_0)$ 
	in \eqref{eq:def_OLNE_ply_2}\,.
\end{definition}

\smallskip

\begin{remark}
	\label{rem:OL_vs_F}
	{\rm Apart from the structure of the underlying control law, it is worth observing that the distinction between F-NE and OL-NE consists also in the characterization of the subset of alternative control laws against which a \emph{candidate} NE should be compared (see the definitions of the feasible sets in Definitions \ref{def:FNE} and \ref{def:OLNE}, respectively). Furthermore, it is worth observing that such admissibility conditions in Definitions \ref{def:FNE} and \ref{def:OLNE}, namely the requirements following the inequalities \eqref{eq:def_FNE} and \eqref{eq:def_OLNE}, respectively, account for the unbounded horizon in the cost functional \eqref{eq:cost_J_i}. As such, they are specific to 
			the infinite-horizon setting and do not possess any counterpart in the finite-horizon formulation. Similar conditions appear in \cite[Sec. 7.4]{Engwerda2005} although in the case of continuous-time LQ games.}
	\hfill $\triangle$
\end{remark}

The objective of this manuscript is to provide a comprehensive analysis of the conditions characterizing both F-NE and OL-NE of the class of dynamic games described by \eqref{eq:lti} and \eqref{eq:cost_J_i}. 
While the study of the two concepts of solution may be carried out in a somewhat parallel way, throughout the paper several connections are established and explicitly discussed \tm{(see, \emph{e.g.} Remarks \ref{rem:notARE} and \ref{rem:comaprison} and Section \ref{sec:sim}).} 

\section{F-NE for infinite-horizon dynamic games}
\label{sec:F_NE}
Focusing first on Definition \ref{def:FNE}, in this section we discuss the characterization of F-NE solutions. 
Consider the following condition, necessary for the existence of a F-NE, that is assumed to hold throughout this section.

\vspace{.7em}

\begin{assumption}
	The pair $\left(A, \begin{bmatrix} B_1 & B_2 \end{bmatrix} \right)$ is stabilizable\tm{, \emph{i.e.} there exists a matrix $K$ of appropriate dimension such that $\sigma(A+\begin{bmatrix} B_1 & B_2) \end{bmatrix}K) \subset \mathbb{C}_d$}.
	\label{as:stabilizability_F-NE}
\end{assumption}

Assumption \ref{as:stabilizability_F-NE} is a necessary requirement to ensure that the feasible set in Definition \ref{def:FNE} is not empty. The verification of such a condition may be performed via a standard rank condition on the reachability matrix for the combined effect of the control actions, namely by interpreting the columns of the matrices $B_1$ and $B_2$ as those of a common input matrix $B$. The following statement provides a characterization of F-NE strategies for the dynamic game described by the system \eqref{eq:lti} and the cost functionals \eqref{eq:cost_J_i}, $i=1,2$, in terms of certain coupled algebraic equations.

\smallskip

\begin{theorem} \label{thm:F_NE}
	Consider the system \eqref{eq:lti} and the cost functionals \eqref{eq:cost_J_i}, for $i = 1,2$. Then the pair of control laws
	\begin{equation}\label{eq:u_F_NE}
		u_i^{F}(x) = K_i x\,,
	\end{equation}
	\noindent for $i = 1,2$, constitutes a F-NE if and only if $K_i \in \mathbb{R}^{m_i \times n}$ are such that
	\begin{equation} \label{eq:stability} 
		\sigma(A+B_1 K_1 + B_2 K_2) \subset \mathbb{C}_{d}\,. 
	\end{equation}
	\noindent and solve, together with some $P_i \in \mathbb{R}^{n \times n}$, $P_i = P_i^{\top} \succeq 0$, $i = 1,2$, the matrix equations
	\begin{subequations}\label{eq:ACE}
		\begin{gather}\label{eq:ARE_1}
			\begin{split}
				P_1 =Q_1 + K_1^\top R_1 K_1   +  (A+B_1K_1 + B_2K_2)^\top P_1 (A+B_1K_1 + B_2K_2)  \,, 
			\end{split} 
		\end{gather} 
		\vspace{-20pt} 
		\begin{gather}\label{eq:ARE_2}
			\begin{split} 
				P_2 = Q_2 + K_2^\top R_2 K_{2}  + (A+B_1K_1 + B_2K_2)^\top P_2 (A+B_1K_1 + B_2K_2)\,,
			\end{split}
		\end{gather} 
	\end{subequations}
	\noindent and
	\begin{subequations}\label{eq:K_M}
		\begin{gather}\label{eq:Ki_s}
			M \left[\begin{array}{c} K_1 \\[2pt] K_2 \end{array}\right] = -\left[\begin{array}{c} B_1^{\top}P_1 \\[2pt] B_2^{\top}P_2 \end{array} \right]A\,, 
		\end{gather}
		\noindent with
		\begin{gather}
			M = \left[\begin{array}{cc} R_1 + B_1^{\top}P_1 B_1 & B_1^{\top}P_1 B_2 \\[2pt] B_2^{\top}P_2 B_1 & R_2 + B_2^{\top}P_2 B_2 \end{array} \right]\,.
		\end{gather}
	\end{subequations}
	Furthermore,  the F-NE is such that   $J_i(x_0, u_1^{F},u_2^{F}) = \frac{1}{2} x_0^\top P_i x_0$,
	for $i=1,2$. 
	\hfill $\circ$
\end{theorem}

\smallskip

\textbf{Proof:} 
Consider first the sufficient implication. Suppose that player 2 adheres to the strategy $u_2 = u^F_2(x)$ and let $A_{cl,2} = A + B_2K_2$ with $K_2$ defined in \eqref{eq:ACE}, \eqref{eq:K_M}. \hl{Note that the first row of \eqref{eq:K_M} is 
\begin{equation}\label{eq:first_row_M}
    (R_1 + B_1^{\top}P_1 B_1) K_1 + B_1^{\top}P_1 B_2 K_2 = - B_1^{\top}P_1 A\,,
\end{equation}
\noindent which yields
\begin{equation}\label{eq:K1_explicit}
    K_1 = -(R_1 + B_1^\top P_1 B_1 )^{-1} B_1^\top P_1 A_{cl,2}\,,
\end{equation}
\noindent where invertibility of the matrix $(R_1 + B_1^{\top} P_1 B_1)$ is ensured by positive definiteness of $R_1$ and positive semi-definiteness of $P_1$.}
Then, with the action of player 2 fixed, \hl{the computation of} the best response of player~1 constitutes a (single-player) optimal control problem  described by the dynamics $x(k+1) = (A+B_2 K_2)x(k) + B_1 u_1(k)$ and the cost functional \eqref{eq:cost_J_i} with $i = 1$. \hl{By replacing \eqref{eq:K1_explicit} into \eqref{eq:ARE_1}, and suitably arranging the terms, it follows that the discrete-time ARE (see, e.g., \cite{kuvcera1972discrete})
\begin{equation}\label{eq:DP1} \begin{split}
		P_{1}^{oc} = & A_{cl,2}^\top P_{1}^{oc} A_{cl,2}  + Q_1 - A_{cl,2}^\top P_{1}^{oc}B_1(R_1 + B_1^\top P_{1}^{oc} B_1)^{-1} B_1^\top P_{1}^{oc} A_{cl,2} \,,
	\end{split} 
\end{equation}
\noindent associated with 
the single-player optimal control problem above admits 
a stabilizing \tm{solution} 
$P_1^{oc} = P_1$. Thus, considering the asymptotic stability requirement of the underlying optimal control problem, one has that the optimal feedback is 
}
\begin{equation}\label{eq:u_oc}
	u_1^{oc}(x) = - (R_1 + B_1^\top P_1^{oc} B_1)^{-1}B_1^\top P_1^{oc} A_{cl,2}x = K_1^{oc} x\,.
\end{equation}
\noindent \hl{A direct inspection of \eqref{eq:u_oc} and \eqref{eq:K1_explicit}}
reveals that $u_1^{oc}(x) = u_1^F(x)$ and \eqref{eq:def_FNE_ply_1} follows. Moreover, it follows from standard Dynamic Programming arguments \hl{that $V_1(x) =\frac{1}{2} x^\top P_1^{oc} x$ constitutes the value function corresponding to the optimal control problem obtained when the strategy of the second player is fixed at $u_2= u_2^F(x)$ (see equation (5) of \cite{kuvcera1972discrete}) and, thus, it follows (by definition of the value function) that $J_1(x_0, u_1^F, u_2^F) = \frac{1}{2}x_0^\top P_1 x_0$, since $P_1 = P_1^{oc}$.}
Consider now the converse statement. 
Assume that the pair of linear strategies $(K_1^\star x,K_2^\star x)$ is a F-NE of the game \eqref{eq:lti}, \eqref{eq:cost_J_i}, $i=1,2$. By Definition~\ref{def:FNE}, \eqref{eq:def_FNE_ply_1} holds for all $K_2$ such that $\sigma(A+B_1 K_1^\star + B_2 K_2) \subset \mathbb{C}_{d}$. Hence, for fixed $K_2$, the feedback strategy $K_1^\star x$ is the unique stabilizing optimal control strategy \eqref{eq:u_oc} with $P_1^{oc}$ satisfying \eqref{eq:DP1} \hl{(see Theorems 2 and 6 of \cite{kuvcera1972discrete} and the comment before Theorem 7 \tm{of the same reference})}, which implies that there exists $P_1 = P_1^{oc}$ such that \eqref{eq:ARE_1} holds with $K_1=K_1^\star.$
The proof is concluded via a symmetric reasoning for player 2.
\hfill $\square$

Theorem~\ref{thm:F_NE} characterizes F-NE strategies via the set of \emph{asymptotically stabilizing} solutions of the coupled algebraic equations \eqref{eq:ACE}, \eqref{eq:K_M}. In the following statement we provide sufficient conditions under which a solution to \eqref{eq:ACE}, \eqref{eq:K_M} is indeed asymptotically stabilizing. 

\smallskip

\begin{corollary}
	\label{co:stab_sol_ARE}
	Let $K_1$, $K_2$, $P_1 \succeq 0$, $P_2 \succeq 0$ be a solution to \eqref{eq:ACE}, \eqref{eq:K_M}. If the pair $\left( A, Q_1+Q_2 \right)$ is detectable\footnote{\hl{The property of \emph{detectability} for a dynamical system is intuitively related to the possibility of constructing asymptotically convergent observers for the state of the system, by assigning all the eigenvalues of the error dynamics within the unit disk. More precisely, a pair $(A,C)$ is detectable if and only if
    \begin{equation}\label{eq:PBH_detect}
        {\rm rank} \left( \left[\begin{array}{c} A - \lambda_i I \\ C \end{array} \right] \right) = n, \hspace{1cm} \forall \lambda_i \in \sigma(A): |\lambda_i| \geq 1\,.
    \end{equation}
    {\hlb The above condition is the so-called Popov-Belevitch-Hautus (PBH) test.}
    \noindent
    Equivalently, $(A,C)$ is detectable if and only if for all $\lambda_i \in \sigma(A)$ and $\lambda_i \notin \mathbb{C}_d$ there exists no $v_i \neq 0$ such that $A v_i - \lambda_i v_i = 0$ together with $C v_i = 0$ (or, alternatively, any right eigenvector of $A$ associated to $\lambda_i \notin \mathbb{C}_d$ is such that $C v_i \neq 0$).}}, then $K_1$, $K_2$ are such that \eqref{eq:stability} holds, \emph{i.e.} the strategies $(K_1x, K_2x)$ constitute a F-NE of the game \eqref{eq:lti}, \eqref{eq:cost_J_i}, $i=1,2$.
	\hfill $\circ$
\end{corollary}

\textbf{Proof:} 
The claim is proved by contradiction. Adding \eqref{eq:ARE_1} and \eqref{eq:ARE_2} gives
\begin{multline}
	(A+B_1 K_1 + B_2 K_2)^\top \tilde P (A+B_1 K_1 + B_2 K_2) -  \tilde P + \tilde Q + \sum_{i=1}^2 K_i^\top R_i K_i = 0,
	\label{eq:summed_ARE}
\end{multline}
where $\tilde P = P_1 + P_2$ and $\tilde Q = Q_1 + Q_2$. Assume \hl{then, by contradiction, that} there exists a solution $K_1$, $K_2$, $P_2 \succeq 0$, $P_1 \succeq 0$ of \eqref{eq:ACE}, \eqref{eq:K_M}, which is such that $\sigma(A+B_1 K_1 + B_2 K_2)\not\subset \mathbb{C}_d$. Let $\lambda\notin \mathbb{C}_d$, with corresponding eigenvector $v\neq0$, be an eigenvalue of the closed-loop equilibrium system that is not located in the interior of the unit circle, \hl{whose existence is implied by the absurd hypothesis. Note that, by definition of eigenvalue and eigenvector, it follows that $(A+B_1 K_1 + B_2 K_2)v = \lambda v$ and $v^{\mathsf{H}}(A+B_1 K_1 + B_2 K_2)^{\top} = \lambda^{\mathsf{H}} v^{\mathsf{H}}$, and hence
\begin{equation}\label{eq:combined_left_right_mult}
    v^{\mathsf{H}}(A+B_1 K_1 + B_2 K_2)^{\top}\tilde{P} (A+B_1 K_1 + B_2 K_2)v = \lambda^{\mathsf{H}} v^{\mathsf{H}} \tilde{P} \lambda v = \lambda^{\mathsf{H}} \lambda v^{\mathsf{H}} \tilde{P} v
\end{equation}}  
\noindent Then, pre and post multiplying \eqref{eq:summed_ARE} by $v^\mathsf{H}$ and $v$, respectively, \hl{and using \eqref{eq:combined_left_right_mult}} give 
\begin{equation}
	\left( \lambda ^\mathsf{H} \lambda - 1 \right) v^\mathsf{H} \tilde P v + v^\mathsf{H} \tilde Q v + v^\mathsf{H} \left(\sum_{i=1}^2 K_i^\top R_i K_i \right) v = 0.
	\label{eq:ARE_stability}
\end{equation}
{\hlb Since by assumption $\left( \lambda^\mathsf{H} \lambda - 1 \right) \hl{= (|\lambda|^2-1)} \geq 0$, $\tilde Q \succeq 0$, $\tilde P \succeq 0$ and $K_i^\top R_i K_i \succeq 0$, for $i = 1,2$, the left-hand side of \eqref{eq:ARE_stability} consists of the sum of non-negative terms. Therefore, since the overall sum is equal to zero by the right-hand side of \eqref{eq:ARE_stability}, it follows that each individual term must be equal to zero, hence $v^\mathsf{H} \tilde Q v = 0$ and $v^\mathsf{H} K_i^\top R_i K_i  v = 0$, $i = 1,2$. The latter identities, being squared (semi-)norms equal to zero,} imply $K_i v = 0$, \hl{for $i = 1,2$}, and $\tilde Q v= 0$. {\hlb Furthermore}, the right eigenvector $v$ of $(A+B_1 K_1 + B_2 K_2)$ associated to $\lambda \notin \mathbb{C}_d$ is such that $(A+B_1 K_1 + B_2 K_2)v = A v =  \lambda v$, hence it is {\hlb also} an eigenvalue also of $A$ {\hlb with the property that} $(Q_1 + Q_2)v = 0$. The latter observation contradicts the hypothesis that the pair $\left( A, Q_1+Q_2 \right)$ is detectable, {\hlb since there exists, as a consequence of the absurd hypothesis, an unobservable eigenvalue (since the PBH observability test is not verified) not belonging to $\mathbb{C}_d$}.
\hfill $\square$

\smallskip

\begin{remark}\label{rem:notARE}
	{\rm The equations \eqref{eq:ACE}, \eqref{eq:K_M} appear reminiscent of the classic AREs arising in (single-player) optimal control problems in discrete-time, with additional terms to account for the presence of the opponent. Note, however, that  the matrix algebraic equations arising in Theorem~\ref{thm:F_NE} are \emph{not quadratic} in the unknown variables. This feature is appreciated in \eqref{eq:ACE}, \eqref{eq:K_M} with respect to the unknowns $K_i$ and $P_i$, $i = 1,2$, and it becomes even more evident whenever $K_i$, $i = 1,2$, can be obtained from \eqref{eq:K_M} and replaced in \eqref{eq:ACE}, \emph{i.e.} provided the matrix $M$ is invertible, thus yielding cubic terms in the variables $P_i$. As such, \eqref{eq:ACE}, \eqref{eq:K_M} \emph{are not Riccati equations}. This is in stark contrast not only  to discrete-time optimal control problems (for which Dynamic Programming leads to quadratic matrix equations), but also to (continuous-time) LQ \emph{differential} games (for which F-NE solutions are characterised by coupled quadratic matrix equations). Interestingly, as shown in Theorem \ref{thm:OL_feedback_synth} below, the computation of OL-NE for discrete-time infinite-horizon games revolves around the solution of certain (asymmetric) Riccati equations.}~\hfill $\triangle$
\end{remark}

\section{OL-NE for infinite-horizon dynamic games}
\label{sec:OL_NE}
We now turn our attention to Definition \ref{def:OLNE} and the class of OL-NE solutions for infinite-horizon LQ dynamic games. The analysis of such solutions 
	is performed according to the following structure. First -- by recognizing that at the core of the computation of a Nash equilibrium one finds a standard (single-player) optimal control problem when the opponent has \emph{fixed} a certain (feed-forward) control strategy -- we develop an infinite-horizon extension of the PMP  for discrete-time linear systems, driven by a known exogenous input in Section~\ref{sec:OC_infinite_hor}. This objective is in turn obtained in two steps: first we recall the optimal feedback obtained via Dynamic Programming arguments (Lemma~\ref{lemma:OC_DP}); then we show that such an optimal control law can be obtained, for a given initial condition of the state, as a specific trajectory of the Hamiltonian dynamics arising in PMP (Lemma~\ref{lemma:OC_PMP_equiv_DP}), thus yielding \emph{sufficient conditions} for optimality on the latter as an open-loop strategy. Finally, in Section \ref{sec:pmp}, such boundary conditions that allow to \emph{identify} the optimal trajectory of the Hamiltonian dynamics are translated into the context of dynamic games (Theorem~\ref{thm:OL_NE}), thus transforming PMP into a set of sufficient conditions for OL-NE solutions. 

\subsection{Infinite-horizon affine-quadratic optimal control}
\label{sec:OC_infinite_hor}
Consider a (single-player) optimal control problem described by the dynamics
\begin{equation}
	\label{eq:OC_sys}
	x(k+1) = A x(k) + B u(k) + w(k)\,, \hspace{0.4cm} x(0) = x_0\,,
\end{equation}
subject to the initial condition $x_0 \in \mathbb{R}^{n}$, where $x : \tm{\mathbb{N}_{\ge0}} \rightarrow \mathbb{R}^n$ denotes the state of the system, $u : \tm{\mathbb{N}_{\ge0}} \rightarrow \mathbb{R}^m$ denotes the control input and $w(\cdot) \in \ell_2$ denotes a \emph{known} signal, and the cost functional 
\begin{equation}
	\label{eq:OC_cost}
	J^{oc}(x_0, u(\cdot)) = \frac{1}{2} \displaystyle \sum_{k = 0}^{\infty} x(k)^{\top}Q x(k) + u(k)^{\top} R u(k)\,,
\end{equation}
with $Q = Q^{\top} \succeq 0$ and $R = R^{\top} \succ 0$.
The characterization of the \emph{open-loop} solution of the problem described by \eqref{eq:OC_sys}, \eqref{eq:OC_cost} is obtained below in two steps: first the \emph{feedback, time-varying} optimal control law is constructed via Dynamic Programming arguments; then \emph{boundary conditions} are obtained under which the induced time history of the optimal control can be equivalently generated by a specific trajectory of the underlying Hamiltonian dynamics. While the problem hinted at in the first step has already been tackled in \cite{engwerda_infinityHor}, the novel constructions of the second step allow to obtain an open-loop implementation of the former result, by establishing 
	a connection with a specific trajectory of the corresponding Hamiltonian dynamics. The boundary conditions associated with this specific trajectory will be instrumental as we turn our attention back to the dynamic game setting and OL-NE solutions in Section \ref{sec:pmp}. Furthermore, since an alternative proof is provided for the first step described above, a formal statement (\emph{i.e.} Lemma \ref{lemma:OC_DP}) is dedicated also to such results, which are akin to those in \cite{engwerda_infinityHor}.

\smallskip

\begin{assumption}
	\label{assu:OC_structural}
	Consider \eqref{eq:OC_sys} and \eqref{eq:OC_cost}. Suppose that \emph{(i)}
	the pair $(A,B)$ is reachable\footnote{\hl {The pair $(A,B)$ is reachable if and only if
\begin{equation}\label{eq:PBH_reach}
    {\rm rank} \left( \left[\begin{array}{cc} A - \lambda_i I & B \end{array} \right] \right) = n, \hspace{1cm} \forall \lambda_i \in \sigma(A).
\end{equation}
The property of reachability for the pair $(A,B)$ ensures that the spectrum of the (closed-loop) matrix $A+BK$ can be \emph{arbitrarily} assigned via the selection of the matrix $K$.}}; \emph{(ii)} the pair $(A,Q)$ is detectable; \emph{(iii)} the matrix $A$ is non-singular; \emph{(iv)} given $x_0$ and $w(\cdot)$, there exists at least one control law $\hat{u}(\cdot)$ such that $J^{oc}(x_0,\hat{u}(\cdot))$ is finite.
	\hfill $\circ$
\end{assumption}

\smallskip

\begin{lemma}
	\label{lemma:OC_DP}
	Consider the system \eqref{eq:OC_sys} and the cost functional \eqref{eq:OC_cost}.	Suppose that Assumption \ref{assu:OC_structural} holds.
	Let $P \in \mathbb{R}^{n \times n}$, $P = P^{\top} \succ 0$ be the stabilizing solution to the \emph{discrete-time Riccati equation}
	\begin{equation}
		\label{eq:OC_Riccati}
		P = Q + A^{\top}P A - A^{\top}P B (R+B^{\top}PB)^{-1} B^{\top}P A\,, 
	\end{equation}
	and let $b(\cdot)$ be defined as
	\begin{equation}
		\label{eq:OC_b_def}
		b(k) = \displaystyle \sum_{h = k}^{\infty} (A_{\rm cl}^{\top})^{h-k+1} P w(h)\,,
	\end{equation}
	with $A_{\rm cl}^{\top} = A^{\top}(I-B(R+B^{\top}PB)^{-1}B^{\top}P)^{\top}$.
	Then the optimal control law is given by 
	\begin{equation}
		\label{eq:OC_u_star}
		\begin{split}
			u^{\star}(k) = -(R+B^{\top}PB)^{-1}B^{\top}P A x(k) 
			- (R+B^{\top}P B)^{-1} B^{\top}(P w(k) + b(k+1))\,.
		\end{split}
	\end{equation}
	Moreover the minimum cost is given by
	\begin{equation}
		\label{eq:OC_optimal_cost}
		J^{oc}(x_0, u^{\star}(\cdot)) = \frac{1}{2}x_0^{\top}P x_0 + x_0^{\top}b(0) + \hl{c_0}\,,
	\end{equation}
	\hl{with $c_0 = \frac{1}{2}\sum_{k=0}^{\infty} w(k)^{\top}P w(k) + 2 w(k)^{\top}b(k+1) - (B^{\top}Pw(k)+B^{\top}b(k+1))^{\top}(R+B^{\top}PB)^{-1}(B^{\top}Pw(k)+B^{\top}b(k+1))$.}
	\hfill $\circ$
\end{lemma}

\smallskip

\textbf{Proof:}
The claim follows from Dynamic Programming arguments. Let 
\begin{equation}
	\label{eq:OC_candidate_value_fct}
	V(k,x(k)) = \frac{1}{2}x(k)^{\top} P x(k) + x(k)^{\top}b(k) + c(k)
\end{equation}
be a \emph{candidate value function} and consider the corresponding Bellman equation
\begin{equation}
	\label{eq:OC_bellman_eq}
	\begin{split}
		V(k,x(k)) = & \displaystyle \min_{u(k)}\Big\{  \frac{1}{2}x(k)^{\top}Q x(k) + \frac{1}{2}u(k)^{\top}R u(k) \\[3pt] & + V(k+1,A x(k) + B u(k) + w(k)) \Big\}\,,
	\end{split}
\end{equation}
which must be verified for all $(k,x(k)) \in\tm{\mathbb{N}_{\ge0}}\times \mathbb{R}^{n}$ {\hlb (see, e.g., \cite[Sec. 3.2]{bertsekas2012dynamic})}.
Straightforward computations show that $u^{\star}(k)$ defined in \eqref{eq:OC_u_star} is the argument that attains the minimum of the right-hand side of \eqref{eq:OC_bellman_eq}.
Replacing  $u^{\star}$ into the latter equation then leads (via straightforward although tedious derivations) to a non-homogeneous quadratic, time-varying function of $x(k)$.
While the (matrix) coefficient of the quadratic terms is zeroed provided the discrete-time Riccati \eqref{eq:OC_Riccati} holds, the coefficient of the linear terms is equal to zero provided $b(\cdot)$ satisfies the difference equation
\begin{equation}
	\label{eq:OC_b_k_eq}
	b(k) = A_{\rm cl}^{\top}b(k+1) + A_{\rm cl}^{\top}P w(k)\,.
\end{equation}


Item (iv) of Assumption~\ref{assu:OC_structural} implies that
$\lim_{k \rightarrow \infty}u^{\star}(k) = 0$, by positive definiteness of $R$. \hl{Moreover, \tm{
item (iv) also implies 
that $\lim_{k \rightarrow \infty} Q x(k) = 0$, which together with item (ii) of Assumption~\ref{assu:OC_structural} implies} that $\lim_{k \rightarrow \infty} x(k) = 0$\tm{.  
Note that  only such a trajectory} is 
compatible with $Q x\tm{(k)}$ tending to zero due to the \tm{aforementioned} detectability assumption}, {\hlb since the latter ensures that all the unobservable modes are asymptotically stable while the observable modes (via the output matrix $Q$) are guaranteed to converge to the origin by the fact that $Q x(k)$ tends to zero as time diverges}. As a consequence, by inspecting the structure of $u^{\star}(k)$ in \eqref{eq:OC_u_star}, one has that $\lim_{k \rightarrow \infty}B^{\top}b(k) = 0$, hence, since $w(k)$ tends to zero 
by definition,
there exists a function $\pi : \tm{\mathbb{N}_{\ge0}} \rightarrow \mathbb{R}^{m}$ with the properties that $B^{\top}b(k) = \pi(k)$ for all $k$ and $\lim_{k \rightarrow \infty} \pi(k) = 0$.
Therefore, by rewriting \eqref{eq:OC_b_k_eq} as $b(k+1) = A_{\rm cl}^{-\top}b(k) - Pw(k)$ and iterating on the reasoning above for $n$ steps, it follows that 
\begin{equation*}
		\begin{split}
			\left[\begin{array}{c} 0 \\ 0 \\ \vdots \\ 0 \end{array}\right] & =   \left[\begin{array}{cccc} 0 & 0 & \ldots & 0 \\-B^{\top}P & 0 & \ldots & 0 \\ \vdots & \ddots &  & \vdots \\-B^{\top}(A_{\rm cl}^{-\top})^{n-2}P & -B^{\top}(A_{\rm cl}^{-\top})^{n-3}P & \ldots & -B^{\top}P \end{array} \right]\left[\begin{array}{c} w(k) \\ w(k+1) \\ \vdots \\ w(k+n-2) \end{array} \right] \\& \qquad +\left[\begin{array}{c} B^{\top} \\ B^{\top}(A_{\rm cl}^{-\top}) \\ \vdots \\B^{\top}(A_{\rm cl}^{-\top})^{n-1} \end{array}\right] b(k) - \left[\begin{array}{c} \pi(k) \\ \pi(k+1) \\ \vdots \\ \pi(k+n-1) \end{array} \right]
			\\[5pt] & =: \mathcal{O}(A_{\rm cl}^{-\top},B^{\top})b(k) + \mathcal{M} W(k) - \Pi(k)\,,
		\end{split}
	\end{equation*}
holds for any $k$. 
By the Woodbury matrix identity\footnote{Note that since $R\succ 0$ and $P\succeq 0$, $(R+B^\top PB) \succ 0$. Thus, the inverse in \eqref{eq:OC_identity} exists, since both the identity matrix and $(R+B^\top PB)$ are invertible (see e.g. \cite[Chapter 2]{golub}).} 
\begin{equation}
	\label{eq:OC_identity}
	(I+B R^{-1} B^{\top}P )^{-1} = I - B(R+B^{\top}P B)^{-1} B^{\top} P\,. 
\end{equation}
Thus, it follows from item (iii) of Assumption \ref{assu:OC_structural}, \hl{ensuring invertibility of the matrix $A$,} and the structure of $A_{cl}$ that $A_{cl}^{-1}$ exists. \hl{In fact, by recalling the expression of $A_{cl}$ in the statement and by using the identity \eqref{eq:OC_identity}, it follows that $A_{cl}$ can be equivalently written as $A_{cl} = (I-B(R+B^{\top}PB)^{-1}B^{\top}P)A = (I + B R^{-1}B^{\top}P)^{-1}A$, the inverse of which is immediately computed as $A_{cl}^{-1} = A^{-1}(I + B R^{-1}B^{\top}P)$.}  Moreover, note that item (i) of Assumption \ref{assu:OC_structural}  implies reachability of the closed-loop (feedback) system $(A_{\rm cl},B)$ and this, in turn, implies reachability of the pair $(A_{\rm cl}^{-1},B)$. \hl{In fact, suppose that the pair $(A_{cl},B)$ is reachable, hence ${\rm rank}([A_{cl} - \lambda_i I \hspace{0.2cm} B]) = n$ for all $\lambda_i \in \sigma(A_{cl})$. Thus
\begin{equation}\label{eq:PBH_for_inverse}
\begin{split}
    n = & {\rm rank}\left(\left[\begin{array}{cc} A_{cl} - \lambda_i I & B \end{array} \right] \right) = {\rm rank}\left( \left[\begin{array}{cc} A_{cl} - \lambda I & B \end{array} \right] \left[\begin{array}{cc} \lambda_i^{-1} A_{cl}^{-1} & 0 \\ 0 & I \end{array}\right] \right) \\[2pt] = & {\rm rank}\left(\left[\begin{array}{cc} \dfrac{1}{\lambda_i} I - A_{cl}^{-1} & B \end{array} \right] \right)\,,
    \end{split}
\end{equation}
\noindent where the multiplication on the second identity does not modify the rank and the last rank is equivalent to the reachability property of the pair $(A_{cl}^{-1},B)$ since if $\lambda_i \in \sigma(A_{cl})$ then $1/\lambda_i \in \sigma(A_{cl}^{-1})$.} By duality, the latter property ensures observability\footnote{\tm{The pair $(A,C)$ is observable if and only if
$
    {\rm rank} \left( \left[\begin{array}{c} A - \lambda_i I \\ C \end{array} \right] \right) = n, \hspace{1cm} \forall \lambda_i \in \sigma(A).
$}} 
of the pair $(A_{\rm cl}^{-\top},B^{\top})$, and hence the matrix $\mathcal{O}(A_{\rm cl}^{-\top},B^{\top})$ is full column rank.
Therefore, $b(k)$ can be uniquely expressed as $b(k) = \mathcal{O}(A_{\rm cl}^{-\top},B^{\top})^{\dagger}(\Pi(k)-\mathcal{M} W(k))$,
and hence $\lim_{k \rightarrow \infty}b(k) = 0$ follows from the asymptotic properties of $w$ and $\pi$.
Therefore, as a result of the latter boundary condition and by considering \eqref{eq:OC_b_k_eq}, the explicit solution is yielded by \eqref{eq:OC_b_def}, 
by the fact that $w \in \ell_2$.
Finally, the last claim follows by observing that the constant term of \eqref{eq:OC_bellman_eq}, with $u = u^{\star}$, with respect to $x(k)$ is zeroed for all $k$ by the selection of $c(\cdot)$ that satisfies the difference equation $c(k) = c(k+1) + \frac{1}{2} w(k)^{\top}P w(k) + 
w(k)^{\top}b(k+1) - \frac{1}{2}(B^{\top}Pw(k)+B^{\top}b(k+1))^{\top}(R+B^{\top}PB)^{-1}(B^{\top}Pw(k)+B^{\top}b(k+1))$, with the boundary condition $\lim_{k \rightarrow \infty} c(k) = 0$ by definition of value function. \hl{Thus, given the latter boundary condition, the (backward) difference equation governing $c(\cdot)$ immediately yields $c(0) = c_0$ via standard arguments related to the explicit solution to affine discrete-time systems.}
\hfill $\square$

\smallskip

\begin{remark}
	\label{rem:OC_u_star_k}
	{\rm The optimal control law $u^{\star}(\cdot)$ in \eqref{eq:OC_u_star} can be equivalently written as a function at the current time $k$ alone by exploiting the difference equation \eqref{eq:OC_b_k_eq}, namely
		\begin{equation}
			\label{eq:OC_u_star_k_alone}
			\begin{split}
				u^{\star}(k) =  -(R+B^{\top}PB)^{-1}B^{\top}P A x(k) - (R+B^{\top}P B)^{-1} B^{\top}A_{\rm cl}^{-\top} b(k)\,.
			\end{split}
		\end{equation}
		It is worth observing that the \emph{anticipative} nature of the control law with respect to the exogenous signal $w$ is preserved since the computation of $b(k)$ via \eqref{eq:OC_b_def} 
		\tm{requires 
  knowledge} of $w(h)$ for all $h \geq k$.}
	\hfill $\triangle$
\end{remark}

The objective of the following result is to relate the time history induced by the optimal control law \eqref{eq:OC_u_star} (computed via Dynamic Programming arguments) with a specific trajectory of the underlying Hamiltonian dynamics (arising when approaching the problem from the perspective of PMP).
To this end, 
define the Hamiltonian function for the optimal control problem \eqref{eq:OC_sys}, \eqref{eq:OC_cost},
\begin{equation}
	\label{eq:OC_hamilt_fct}
	\begin{split}
		&\mathcal{H}(x(k),\lambda(k+1),u(k),k) =  \frac{1}{2}x(k)^{\top}Q x(k) \\ & + \frac{1}{2}u(k)^{\top} R u(k) 
		+ \lambda(k+1)^{\top}(A x(k) + B u(k) + w(k))\,.
	\end{split}
\end{equation}
Utilizing PMP arguments results, \hl{for all $k \in \mathbb{N}_{\geq 0}$,} in the 
Hamiltonian dynamics
\begin{subequations}
	\label{eq:OC_hamilt_dyn}
	\begin{align}
		x(k+1) &  = Ax(k) - BR^{-1} B^{\top} \lambda(k+1) + w(k)\,, \label{eq:OC_hamilt_dyn_x} \\[3pt]
		\lambda(k) & = Q x(k) + A^{\top} \lambda(k+1)\,. \label{eq:OC_hamilt_dyn_lam} 
	\end{align}
\end{subequations}

\smallskip

\begin{lemma}
	\label{lemma:OC_PMP_equiv_DP}
	Consider the system \eqref{eq:OC_sys} together with the cost functional \eqref{eq:OC_cost}.
	Suppose that the hypotheses of Lemma \ref{lemma:OC_DP} hold.
	Let $P \in \mathbb{R}^{n \times n}$, $P = P^{\top} \succeq 0$ be the stabilizing solution to \eqref{eq:OC_Riccati} and let $b(\cdot)$ be determined by \eqref{eq:OC_b_def}, and hence such that $\lim_{k \rightarrow \infty}b(k) = 0$.
Then the set defined by \hl{$\{ (\lambda(k), x(k)) : \lambda(k) = P x(k) + b(k), k = 0,1,2, ... \}$} is invariant for \eqref{eq:OC_hamilt_dyn}, for all $k\geq 0$. Moreover,  the optimal control law is given by 
	\begin{equation}
		\label{eq:OC_equivalent_u}
		u^{P}(k) = -R^{-1}B^{\top}\lambda(k+1) = u^{\star}(k)\,,
	\end{equation}
	for all $k \geq 0$, 
	with $u^{\star}$ defined in \eqref{eq:OC_u_star}, where $\lambda(k)$ is the solution of \eqref{eq:OC_hamilt_dyn}, subject to the initial condition $x(0) = x_0$ and $\lambda(0) = Px_0 + b(0)$.
	\hfill $\circ$
\end{lemma}

\smallskip

\textbf{Proof:}
To begin with, \hl{invariance of the identity $\lambda(k) = P x(k) + b(k)$, for all $k \in \mathbb{N}$, with respect to the dynamics \eqref{eq:OC_hamilt_dyn} is shown. {\hlb This is achieved by induction noting that the identity holds at $k = 0$ via a proper initialization of $\lambda(0)$ and $x(0)$, \emph{i.e.} with $x(0)$ and $\lambda(0)$ as specified in the statement. Then, by assuming that $\lambda(k) = P x(k) + b(k)$ holds, it is shown that the same identity is verified also at $k+1$. Towards this end, by item \emph{(iii)} of Assumption \ref{assu:OC_structural}, the difference equation \eqref{eq:OC_hamilt_dyn} may be written in the more standard form of 
\begin{subequations}
	\label{eq:OC_hamilt_dyn_same_time}
	\begin{align}
		x(k+1) &  = (A+B R^{-1} B^{\top}A^{-\top}Q)x(k) - B R^{-1}B^{\top} A^{-\top} \lambda(k) + w(k)\,, \label{eq:OC_hamilt_dyn_x_proof} \\[3pt]
		\lambda(k+1) & = -A^{-\top}Qx(k) + A^{-\top}\lambda(k)\,, \label{eq:OC_hamilt_dyn_lam_proof} 
	\end{align}
\end{subequations}
\noindent together with the \emph{forward time} difference equation for $b(\cdot)$ in \eqref{eq:OC_b_k_eq}, namely $b(k+1) = A_{\rm cl}^{-\top} b(k) - P w(k)$. Then, the first claim is obtained provided one is able to show that
\begin{equation}\label{eq:one_step_invariance}
0 =  - \lambda(k+1) + P x(k+1) + b(k+1)     
\end{equation}
along the trajectories of \eqref{eq:OC_hamilt_dyn_same_time}, whenever $\lambda(k) = P x(k) + b(k)$ is verified. Towards this end, observe that the right-hand side of \eqref{eq:one_step_invariance} becomes
\begin{equation}\label{eq:invariance_proof}
    \begin{split}
 & A^{-\top} Q x(k) - A^{-\top} \lambda(k) + P (A + B R^{-1} B^{\top}A^{-\top}Q) x(k) - P B R^{-1} B^{\top} A^{-\top} \lambda(k) \\[2pt] & \hspace{0.4cm} + P w(k) + A_{\rm cl}^{-\top} b(k) - P w(k) \\[2pt] = \hspace{0.2cm} & (A^{-\top}Q + PA + P B R^{-1}B^{\top}A^{-\top}Q)x(k) - (I + P B R^{-1} B^{\top})A^{-\top} \lambda(k) + A_{\rm cl}^{-\top} b(k) \\[2pt]
  = \hspace{0.2cm} & (A^{-\top}Q + PA + P B R^{-1}B^{\top}A^{-\top}Q)x(k) - A_{\rm cl}^{-\top} \lambda(k) + A_{\rm cl}^{-\top} b(k)
  \\[2pt]
  = \hspace{0.2cm} & (A^{-\top}Q + PA + P B R^{-1}B^{\top}A^{-\top}Q - (I + P B R^{-1} B^{\top})A^{-\top} P) x(k) = 0
    \end{split}
\end{equation}
\noindent where the third equality is obtained by replacing $\lambda(k) = P x(k) + b(k)$, while the last identity follows by observing that the matrix at the last line of \eqref{eq:invariance_proof} is equivalent to \eqref{eq:OC_Riccati}, and hence equal to zero. Such equivalence can be established via lengthy, although straightforward, manipulations revolving around the use of \eqref{eq:OC_identity}.}
The second claim is obtained by first considering the system \eqref{eq:OC_sys} in closed-loop with \eqref{eq:OC_u_star} described (utilizing ~\eqref{eq:OC_identity}) by the equation\footnote{Note that in the proof of Lemma \ref{lemma:OC_PMP_equiv_DP} the superscripts $P$ and $D$ are used to indicate the state trajectory ensuing from \eqref{eq:OC_hamilt_dyn} (PMP) and \eqref{eq:OC_sys}, \eqref{eq:OC_u_star} (DP), respectively.}
\begin{equation}
	\label{eq:OC_closed_loop_DP}
	\begin{split}
		x^{D}(k+1) 
		& = (I+BR^{-1}B^{\top}P)^{-1}(Ax^{D}(k) + w(k)) \\[2pt] & \hspace{0.4cm} - B(R+B^{\top}PB)^{-1}B^{\top}b(k+1)\,.
	\end{split}
\end{equation}
Therefore, exploiting again~\eqref{eq:OC_identity}, it follows that  
\begin{equation*}
	\begin{split}
		(I +B R^{-1}  &B^{\top}P)x^{D}(k+1)  
		\\& = \, A x^{D}(k) + w(k)
		+BR^{-1}B^{\top}b(k+1)\,,
	\end{split}
\end{equation*}
{\hlb Let $x^{P}(\cdot)$ instead denote the solution of \eqref{eq:OC_hamilt_dyn_x} \hl{and suppose that} $\lambda(k+1)$ is replaced by $P x^{P}(k+1) + b(k+1)$. It follows from \eqref{eq:OC_hamilt_dyn_x} that
\begin{equation}
	\label{eq:OC_cl_sys_xP}
	\begin{split}
		x^{P}(k+1) = & \, (I + BR^{-1}B^{\top}P)^{-1}(A x^{P}(k) + w(k)) \\[3pt] & -(I + BR^{-1}B^{\top}P)^{-1}BR^{-1}B^{\top} b(k+1)\,,
	\end{split}
\end{equation}}
\noindent which in turn, by inspecting \eqref{eq:OC_cl_sys_xP} and since $x^{D}(0) = x^{P}(0) = x_0$, implies that $u^{\star}(\cdot)$ and $u^{P}(\cdot)$ induce an identical (optimal) closed-loop trajectory.
\hfill $\square$

	\begin{remark}\label{rem:sufficient_cond}
		The 
		implications of \eqref{eq:OC_equivalent_u} are explicitly discussed here. While conditions derived via Dynamic Programming arguments (such as the ones of Lemma~\ref{lemma:OC_DP} in the context of \eqref{eq:OC_sys}, \eqref{eq:OC_cost}) provide necessary and sufficient conditions for optimality, those based instead on Pontryagin's principle typically yield only necessary requirements. However, the merit of Lemma~\ref{lemma:OC_PMP_equiv_DP} consists in determining certain \emph{boundary conditions} for the Hamiltonian dynamics~\eqref{eq:OC_hamilt_dyn} with the property that the ensuing output trajectory coincides \emph{for any time} with the optimal solution derived in Lemma~\ref{lemma:OC_DP}. As a consequence, the conditions derived via PMP become also sufficient at least for the given initial condition $x_0$ for which the boundary conditions have been determined. More precisely, as entailed by the statement of Lemma~\ref{lemma:OC_PMP_equiv_DP}, such boundary conditions are enforced by determining the costate trajectory that asymptotically converges to the origin, namely such that $\lim_{k \rightarrow \infty} \lambda(k) = 0$. In fact, the latter claim is derived by noting that 
			\begin{equation}\label{eq:lambda_invar_time_varying}
				\lambda(k) = P x(k) + b(k)\,,
			\end{equation}
			\noindent for all $\tm{k  \in \mathbb{N}_{\geq 0}}$, and that $x(k)$ and $b(k)$ tend 
            to zero as $k$ tends to infinity.
		\hfill $\triangle$
\end{remark}

\subsection{OL-NE via Riccati equations}\label{sec:pmp}

The key conclusion of            Lemma~\ref{lemma:OC_PMP_equiv_DP} and of Remark~\ref{rem:sufficient_cond} can be summarized as follows: for fixed $x_0 \in \mathbb{R}^n$, the Hamiltonian dynamics~\eqref{eq:OC_hamilt_dyn} yield sufficient conditions for optimality provided one selects the trajectory of~\eqref{eq:OC_hamilt_dyn} such that $\lim_{k\rightarrow\infty} \lambda(k) = 0$. This observation is instrumental as we
now turn our attention back to dynamic games and OL-NE strategies.
Consider the class of discrete-time systems characterized by a non-singular matrix $A$, and with certain structural properties as prescribed in the following assumptions.

\smallskip

\begin{assumption}
	\label{assu:A_nonsing}
	Consider the system \eqref{eq:lti}.
	The matrix $A$ is non-singular, \emph{i.e.} $\sigma(A) \cap \{0\} = \emptyset$.
	\hfill $\circ$    
\end{assumption}

\smallskip

\begin{assumption}
	\label{assu:structural_games}
	Consider the system \eqref{eq:lti} and the cost functionals \eqref{eq:cost_J_i}. The pairs $(A,B_i)$ and $(A,Q_i)$ are reachable and detectable, respectively, for $i = 1,2$.
	\hfill $\circ$ 
\end{assumption} 

The constructions of this section heavily rely on the properties of the matrix $H \in \mathbb{R}^{3n \times 3n}$ defined, under the hypothesis that Assumption~\ref{assu:A_nonsing} holds, by
\begin{equation}
	\label{eq:def_mat_H}
	H \! := \! \left[\begin{array}{ccc}
		A + S_1A^{-\top}Q_1 + S_2A^{-\top}Q_2 & -S_1A^{-\top} & -S_2A^{-\top} \\
		-A^{-\top}Q_1 & A^{-\top} & 0 \\
		-A^{-\top}Q_2 & 0 & A^{-\top}
	\end{array}\right]\!\!.
\end{equation}
with $S_i = B_i R_i^{-1} B_i^{\top}$, $i = 1,2$. 

\smallskip




\smallskip

\begin{theorem}
	\label{thm:OL_NE}
	Consider the system \eqref{eq:lti} and the cost functionals \eqref{eq:cost_J_i}, for $i = 1, 2$.
	Suppose that Assumptions \ref{assu:A_nonsing} and \ref{assu:structural_games} hold and fix any $x_0 \in \mathbb{R}^{n}$.
	Suppose that the state/costate dynamics
	\begin{equation}
		\label{eq:state_costate_dyn}
		\left[\begin{array}{c} x(k+1) \\ \lambda_1(k+1) \\ \lambda_2(k+1)  \end{array} \right] = H \left[\begin{array}{c} x(k) \\ \lambda_1(k) \\ \lambda_2(k)  \end{array} \right]\,,
	\end{equation}
	\noindent admit a solution, denoted $(x^{OL}, \lambda_1^{OL}, \lambda_2^{OL})$, such that $x(0) = x_0$ and $\lim_{k \rightarrow \infty}\lambda_i(k) = 0$. 
	Then the underlying dynamic game possesses an  OL-NE described by
	\begin{equation}
		\label{eq:u_OL_NE}
		u^{OL}_i(k) = -R_i^{-1}B_i^{\top}\lambda_i^{OL}(k+1)\,,
	\end{equation}
	for $i = 1,2$. 
	\hfill $\circ$
\end{theorem}

\smallskip

\textbf{Proof:}
The claim is shown by borrowing, symmetrically for each player, the constructions of Lemma \ref{lemma:OC_PMP_equiv_DP} with respect to the Hamiltonian functions
\begin{equation}
	\label{eq:hamiltonian}
	\begin{split}
		\mathcal{H}_i = &\,\,\frac{1}{2}x(k)^\top Q_ix(k) + \frac{1}{2}u_{i}(k)^\top R_iu_{i}(k) \\
		&+ \lambda_{i}(k+1)^\top\left(Ax(k) + B_iu_{i}(k) + B_j\phi_j^{\star}(k,x_0)\right),
	\end{split}
\end{equation}
for $i = 1, 2$, for $j = 1, 2$ and with $i \neq j$, namely as in \eqref{eq:OC_hamilt_fct} with the role of $w(\cdot)$ played by $\phi_j^{\star}(\cdot,x_0)$. 
The conclusions of Lemmas \ref{lemma:OC_DP} and \ref{lemma:OC_PMP_equiv_DP}, together with the comments in Remark~\ref{rem:sufficient_cond}  imply that the above underlying (partial) optimal control problem has the unique solution described by \eqref{eq:u_OL_NE} for each player, with the individual costate variables satisfying
\begin{equation}
	\label{eq:costate}
	\lambda_i(k) = Q_ix(k) + A^\top\lambda_i(k+1)\,,
\end{equation}
provided the boundary condition $\lim_{k \rightarrow \infty} \lambda_i(k) = 0$ holds, consistently \tm{with} 
the conclusions of Lemma \ref{lemma:OC_PMP_equiv_DP} in the case of a single-player. The latter is precisely the property that characterizes the process $(x^{OL}, \lambda_1^{OL}, \lambda_2^{OL})$ and hence it holds by assumption.
\tm{By} replacing the control laws \eqref{eq:u_OL_NE} into \eqref{eq:lti} a composite description of the evolution of the individual costate variables is obtained as
\begin{equation}\label{eq:composite_ham_dyn_mixed_times}
	\left[\begin{array}{c}
		x(k+1) \\
		\lambda_1(k) \\
		\lambda_2(k)
	\end{array}\right] = \left[\begin{array}{ccc}
		A & -S_1 & -S_2 \\
		Q_1 & A^\top & 0 \\
		Q_2 & 0 & A^\top
	\end{array}\right]
	\left[\begin{array}{c}
		x(k) \\
		\lambda_1(k+1) \\
		\lambda_2(k+1)
	\end{array}\right].
\end{equation}
By arranging the equations in \eqref{eq:composite_ham_dyn_mixed_times} in order to induce consistent time instants on each side, one immediately obtains
\begin{equation}
	\hspace{-0.3cm}  \label{eq:composite_ham_dyn_pre_inverse}
	\left[\begin{array}{c}
		x(k+1) \\
		\lambda_1(k+1) \\
		\lambda_2(k+1)
	\end{array}\right] \! = \! \left[\begin{array}{ccc}
		I & S_1 & S_2 \\
		0 & -A^\top & 0 \\
		0 & 0 & -A^\top
	\end{array}\right]^{-1} \! \! \left[\begin{array}{ccc}
		A & 0 & 0 \\
		Q_1 & -I & 0 \\
		Q_2 & 0 & -I
	\end{array}\right] \! \!
	\left[\begin{array}{c}
		x(k) \\
		\lambda_1(k) \\
		\lambda_2(k)
	\end{array}\right]\!\!.
\end{equation}
The proof is concluded by noting that, taking the matrix inverse,
the state/costate dynamics \eqref{eq:composite_ham_dyn_mixed_times} yields \eqref{eq:state_costate_dyn}, since
\begin{equation}
	\label{eq:H_alternative}
	H = \left[\begin{array}{ccc}
		I & S_1A^{-\top} & S_2A^{-\top} \\
		0 & -A^{-\top} & 0 \\
		0 & 0 & -A^{-\top}
	\end{array}\right]
	\left[\begin{array}{ccc}
		A & 0 & 0 \\
		Q_1 & -I & 0 \\
		Q_2 & 0 & -I
	\end{array}\right],
\end{equation}
coincides with the matrix defined in \eqref{eq:def_mat_H}. Therefore, since by the considered boundary conditions of the process $(x^{OL}, \lambda_1^{OL}, \lambda_2^{OL})$ the PMP conditions become sufficient, one has that \eqref{eq:u_OL_NE} \emph{is} the optimal solution of the underlying single-player optimal control problem for a fixed strategy of the opponent. Finally, since the same reasoning applies symmetrically to both players, it can be concluded that \eqref{eq:u_OL_NE}, for $i = 1,2$ constitute open-loop Nash equilibrium strategies.  
\hfill $\square$

The following assumption extends to the discrete-time framework conditions inspired by the \emph{strong stabilizability}  of \cite{Engwerda2007} for continuous-time LQ games.

\smallskip

\begin{assumption}
	\label{assu:spectrum_H}
	The matrix $H$ in \eqref{eq:def_mat_H} possesses $n$ eigenvalues, counted with their multiplicity, in $\mathbb{C}_d$ and $2n$ eigenvalues with modulus greater than or equal to one. Moreover, letting $\mathcal{V}^{-}$ denote the $n$-dimensional stable invariant subspace of $H$, the subspaces $\mathcal{V}^{-}$ and 
		$$ {\rm Im}\left(\left[\begin{array}{cc} 0 & 0 \\I_n & 0 \\0 & I_n \end{array} \right]\right) $$
		\noindent are complementary.
	
	\hfill $\circ$    
\end{assumption}

\noindent	Theorem~\ref{thm:OL_NE} provides sufficient conditions for the existence of an OL-NE \emph{for a fixed initial condition}. Therefore, the control law~\eqref{eq:u_OL_NE} does not in general satisfy the underlying Nash inequalities whenever its performance is compared with feedback strategies (see also the numerical simulations \tm{in} 
Section \ref{sec:sim}), namely it does not constitute a F-NE. Nonetheless, it may be still desirable to express~\eqref{eq:u_OL_NE} as a feedback from the current state, \emph{i.e.} in providing the so-called \emph{feedback synthesis}\footnote{See \cite[Sec. 7.4]{Engwerda2005} for more detailed discussions.} of such OL-NE.
This section is concluded by showing that under 
	Assumption \ref{assu:spectrum_H} there exists a unique pair of OL-NE strategies \eqref{eq:u_OL_NE} permitting a feedback synthesis. This solution is
characterized in Theorem~\ref{thm:OL_feedback_synth},  
via the solution of a set of coupled Riccati equations. 

\smallskip

\begin{theorem}
	\label{thm:OL_feedback_synth}
	Consider the system \eqref{eq:lti} and the cost functionals \eqref{eq:cost_J_i}, for $i = 1, 2$.
	Suppose that the hypotheses of Theorem \ref{thm:OL_NE} and Assumption  \ref{assu:spectrum_H} hold and fix any $x_0 \in \mathbb{R}^{n}$.
	Then the underlying dynamic game possesses a unique OL-NE admitting a feedback synthesis, described by
	\begin{equation}
		\label{eq:u_OL_NE_feedback_synth}
		u^{OL}_i(k) = -R_i^{-1}B_i^{\top}A^{-\top}(-Q_i +P_i^p)x(k) =: K_i^{OL} x(k)\,,
	\end{equation}
	for $i = 1,2$, with $P_i^p \in \mathbb{R}^{n \times n}$ solving the coupled Riccati equations
	\begin{equation}
		\label{eq:OL_ARE}
		\begin{split}
			0 = & A^{-\top}Q_1 + P_1^p(A + S_1A^{-\top}Q_1 + S_2A^{-\top}Q_2) -A^{-\top}P_1^p \\ & 
            - P_1^p S_1 A^{-\top}P_1^p - P_1^pS_2A^{-\top}P_2^p\,, \\[4pt]
			0 = & A^{-\top}Q_2 + P_2^p(A + S_1A^{-\top}Q_1 + S_2A^{-\top}Q_2)-A^{-\top}P_2^p \\ & 
            - P_2^pS_1A^{-\top}P_1^p - P_2^pS_2A^{-\top}P_2^p.
		\end{split}
	\end{equation}
	\hfill $\circ$
\end{theorem}

\textbf{Proof:}
Since the costate variables should converge asymptotically to zero, their initial condition can be determined by restricting the dynamics to the stable invariant subspace, whose existence and dimension are ensured by Assumption \ref{assu:spectrum_H}.
One should determine matrices $P_1^p$, $P_2^p$ and $\Lambda$, not necessarily symmetric, such that
\begin{equation}\label{eq:invariance}
	H
	\left[\begin{array}{c}
		I \\
		P_1^p \\
		P_2^p
	\end{array}\right] = \left[\begin{array}{c}
		I \\
		P_1^p \\
		P_2^p
	\end{array}\right]\Lambda,
\end{equation}
with $\Lambda$ such that $\sigma(\Lambda) \subset \mathbb{C}_d$.
Replacing $H$ defined in \eqref{eq:def_mat_H} into \eqref{eq:invariance} and substituting the first row of the resulting system of equations into the bottom two gives  \eqref{eq:OL_ARE}.
Finally, the structure of \eqref{eq:u_OL_NE_feedback_synth} is obtained by restricting the evolution of the costate $\lambda_i(k+1)$ to the set defined by
\begin{equation}\label{eq:lambda_invar_static}
	\lambda_i(k) = P_i^p x(k)
\end{equation}
\noindent as ensured by the \emph{invariant} subspace characterized in \eqref{eq:invariance}, namely such that $\lambda_i(k+1) = -A^{-\top}Q_i x(k) + A^{-\top}P_i^p x(k)$, for $i = 1, 2$.
\hfill $\square$

\begin{remark}\label{rem:feedback_synthesis}
		{\rm The arguments employed in the proof of Theorem \ref{thm:OL_feedback_synth} show that the (somewhat more abstract) boundary conditions of Theorem \ref{thm:OL_NE} can be more effectively replaced by the initial condition $\lambda_i(0) = P_i^p x_0$, with $P_i^p$, $i = 1,2$, solutions of the asymmetric ARE \eqref{eq:OL_ARE}. Furthermore, it is worth relating the two different expressions for the costate variable, namely \eqref{eq:lambda_invar_time_varying} (although in the context of single-player optimal control) and \eqref{eq:lambda_invar_static}: the former is obtained in terms of a time-varying relation in which the underlying matrix $P$ solves conditions derived via Dynamic Programming arguments; the latter is a static relation between state and costate derived by invariance reasonings on the Hamiltonian dynamics arising in PMP. Therefore, although seemingly similar in their expressions, the two conditions are different in their essence (see also similar results in the context of nonlinear continuous-time differential games in \cite{sassano2022analysis}).}
		\hfill $\triangle$ 
\end{remark}

\smallskip

\begin{remark}\label{rem:relax_assumption_4}
	A milder version of Assumption \ref{assu:spectrum_H} may be stated by requiring, instead, that the spectrum of $H$ contains \emph{at least} $n$ eigenvalues with modulus smaller than one. The conclusions of Theorem \ref{thm:OL_feedback_synth} remain valid apart from \emph{uniqueness} of the OL-NE (a continuous-time counterpart of this aspect can be found in the discussions following Theorem 7.11 of \cite{Engwerda2005}).
	\hfill $\triangle$    
\end{remark}

\smallskip

\begin{remark}
	\label{rem:comparison_AbouKandil}
	{\rm A similar formulation of dynamic games has been addressed in \cite{freiling1999discrete}, in which the  focus is on an elegant characterization of the limiting behavior of the finite-horizon Riccati difference equations as the terminal time tends to infinity.
		Differently from \cite{freiling1999discrete}, 
		we provide sufficient conditions that are derived by tackling directly the infinite-horizon problem and this objective is achieved via a preliminary, and relevant \emph{per se}, consideration of PMP applied to a (single-player) optimal control problem with underlying dynamics influenced by a known exogenous input. 
		Moreover, the conditions \eqref{eq:OL_ARE} are more reminiscent of their continuous-time counterpart than the structure of the algebraic equations derived in \cite{freiling1999discrete}.
		This is obtained by relying on the state/costate extended dynamics \eqref{eq:state_costate_dyn}, with the property that the eigenvalues of the closed-loop equilibrium dynamics are a subset of $\sigma(H)$, instead of the set containing the \emph{inverse} of the eigenvalues of the dynamics matrix of the state/costate system used in \cite{freiling1999discrete}. \hl{ 
  \tm{A similar} problem is tackled over an infinite horizon in \cite{kremer}, in which a solution based on a combination of dynamic programming and \tm{an} operator approach is proposed.} }
	\hfill $\triangle$     
\end{remark}

\smallskip

\begin{remark}\label{rem:comaprison}
	{\rm Theorem~\ref{thm:OL_feedback_synth} provides conditions under which there exist constant matrices $K_i^{OL} \in \mathbb{R}^{m_i \times n}$ such that a pair of OL-NE strategies can be implemented as a feedback law, namely $\phi_i^{\star}(k,x_0) = K_i^{OL} x^{\star}(k)$, for $i = 1, 2$ and for all $k \in \mathbb{N}$, where $x^{\star}$ satisfies $x^{\star}(k+1) = (A + B_1 K_1^{OL} + B_2 K_2^{OL})x^{\star}(k)$, $x(0) = x_0$.
		However, recall from Remark~\ref{rem:OL_vs_F} that the strategies \eqref{eq:u_OL_NE_feedback_synth} do not constitute a F-NE solution. 
		Namely, there may exist a gain matrix $\hat{K}_1$ (similarly for the second player) such that
		\begin{equation}
			\label{eq:OL_F_counter_ex}
			J_1(x_0, K_1^{OL}\hat{x}(\cdot), K_2^{OL}\hat{x}(\cdot)) \geq J_1(x_0, \hat{K}_1 \hat{x}(\cdot), K_2^{OL}\hat{x}(\cdot))\,,
		\end{equation}
		where $\hat{x}$ satisfies $\hat{x}(k+1) = (A+B_1 \hat{K}_1 + B_2 K_2^{OL} )\hat{x}(k)$, $\hat{x}(0) = x_0$. This is, in fact, demonstrated in the numerical example considered in Section \ref{sec:sim}.
		Nonetheless, 
		it may still be desirable to compute a feedback synthesis of an OL-NE for several reasons, including a certain margin of robustness guaranteed by the feedback implementation and, more importantly, the possibility of simultaneously computing an OL-NE \emph{for all initial conditions}, without the need to resort to dedicated boundary conditions (see the comments in Remark \ref{rem:feedback_synthesis}).
		It is also worth comparing the coupled matrix equations \eqref{eq:u_OL_NE_feedback_synth}, \eqref{eq:OL_ARE} and \eqref{eq:ACE}, \eqref{eq:K_M} characterizing 
		OL-NE and 
		F-NE solutions, respectively. While the former are quadratic in the (in general asymmetric) variables $P_i^p$, the latter are not quadratic in the symmetric variables $P_i$, $i =1,2$. Hence, depending on the dimension of the system \eqref{eq:lti}, OL-NE strategies admitting a feedback synthesis may be easier to compute and may hence 
		be of practical interest.}
	\hfill $\triangle$     
\end{remark}

\section{Numerical simulations}\label{sec:sim}
In this section the presented results are illustrated 
via a numerical example. 
An example with scalar dynamics is considered to permit an immediate visualization of the results. Namely, consider \eqref{eq:lti} and  \eqref{eq:cost_J_i}, 
$i = 1, 2$, with $A~=~1, B_1~=~2, B_2~=~1, Q_1~=~1, Q_2~=~0.2, R_1~=~2$ and $R_2~=~0.5$.
For the given system and cost parameters there exists a unique F-NE. The corresponding strategies are described by \eqref{eq:u_F_NE} with $(K_1,K_2) = ( -0.3205,-0.1096)$, which satisfy \eqref{eq:ACE}, \eqref{eq:K_M} with $(P_1,P_2) = (1.2854,0.2197)$. Note that the present example is such that Assumption~\ref{assu:spectrum_H} holds. Hence, there exists a unique OL-NE admitting a feedback synthesis. The 
corresponding strategies are given by \eqref{eq:u_OL_NE_feedback_synth} with $(P_1, P_2) = (1.3165, 0.2633)$, $(K_1^{OL},K_2^{OL}) = (-0.3165,-0.1266)$ satisfying \eqref{eq:OL_ARE}.
The time histories of the state, starting from $x(0) =$ 0.95751, corresponding to the F-NE ($u_i=u_i^F$, $i=1,2$, red line) and the OL-NE ($u_i=u_i^{OL}$, $i=1,2$, blue line) are shown in Figure~\ref{fig:x_comparison}. The time histories of the control inputs $u_1(k)$ (solid lines) and $u_2(k)$ (dashed lines) corresponding to the F-NE (red) and OL-NE (blue) strategies are shown in Figure \ref{fig:u_comparison}. 
\tm{Inspecting} Figures~\ref{fig:x_comparison} and \ref{fig:u_comparison} immediately reveals, as expected, that the two set of conditions lead to two different sets of equilibrium strategies, which in turn induce distinct time evolutions of the underlying plant.
\begin{figure}[b!]
	\centering
	\includegraphics[width=0.6\columnwidth]{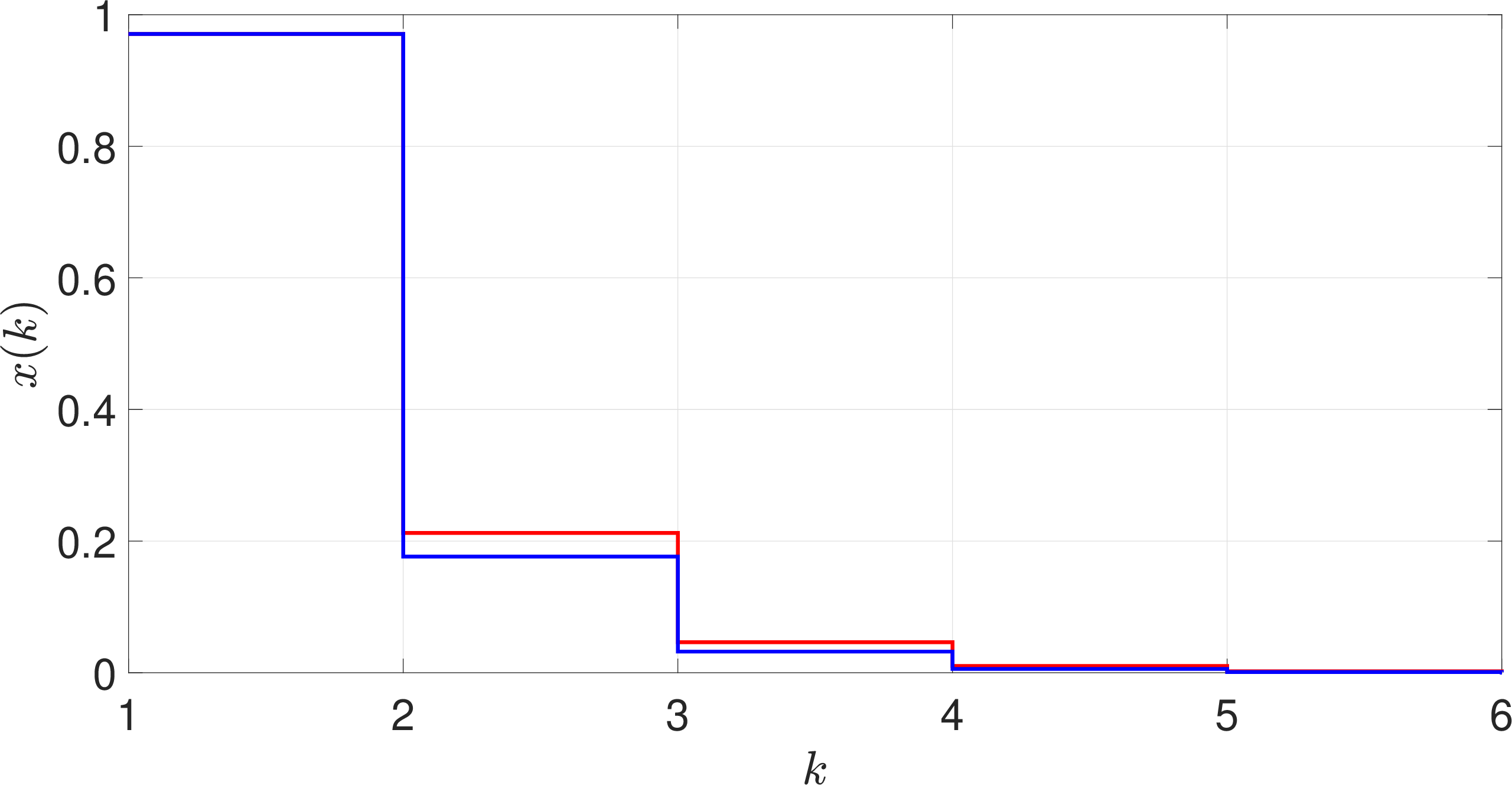}
	\caption{Time histories $x(k)$ corresponding to \eqref{eq:lti} in closed loop with the feedback synthesis of the OL-NE strategies 
		(blue) and with the F-NE strategies
		(red).   
	}
	\label{fig:x_comparison}
\end{figure}
\begin{figure}[htb!]
	\centering
	\includegraphics[width=0.6\columnwidth]{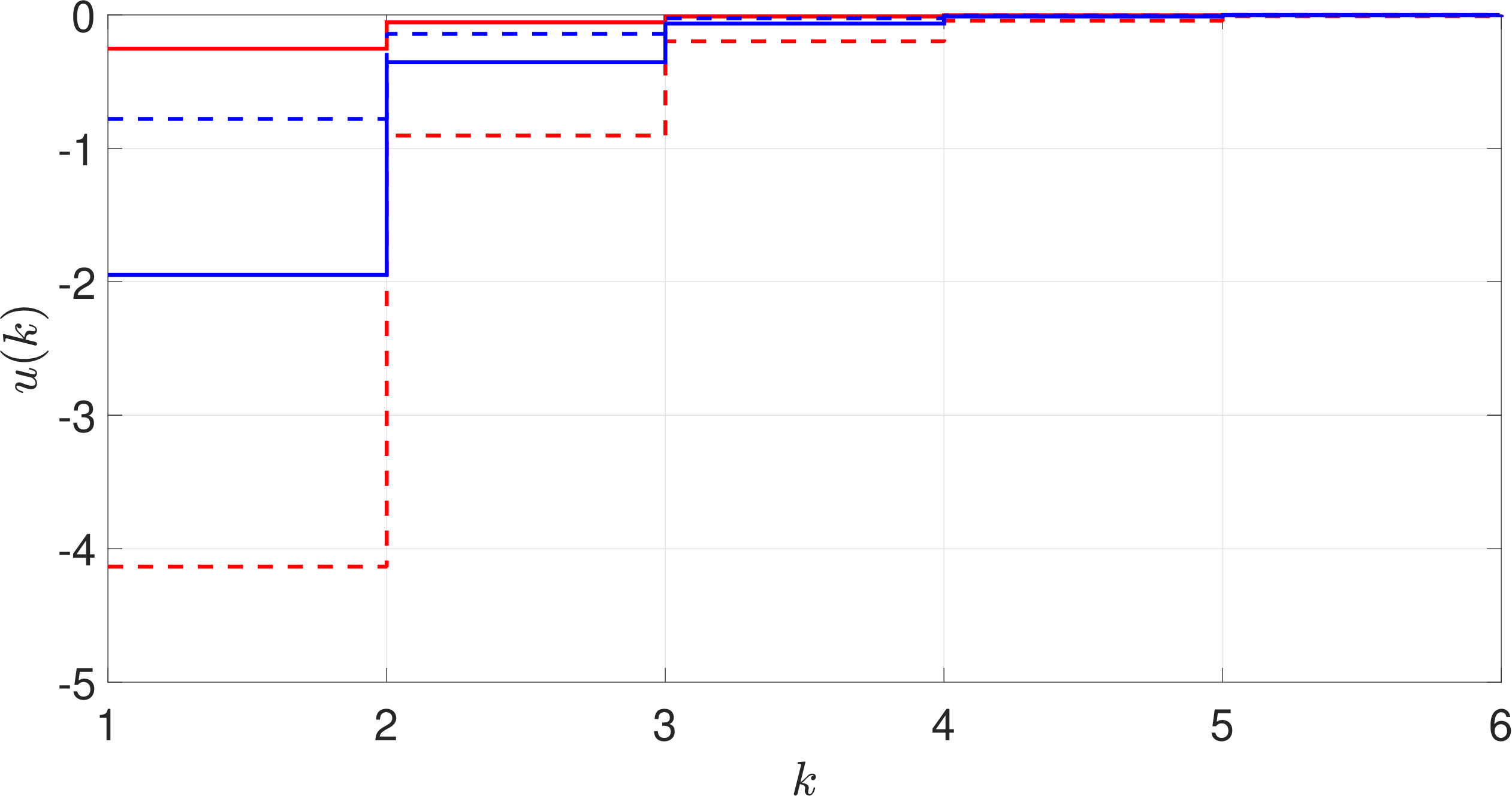}
	\caption{Time histories of $u_1(k)$ (solid lines) and $u_2(k)$ (dashed lines) corresponding to the OL-NE (blue) and the F-NE (red) strategies.}
	\label{fig:u_comparison}
\end{figure}

A different set of numerical simulations is now considered,  with the objective of demonstrating the observations made in Remarks \ref{rem:OL_vs_F} and \ref{rem:comaprison} and graphically illustrating the difference between the Definitions \ref{def:FNE} and \ref{def:OLNE}. To this end, we focus on the pair of OL-NE strategies determined above. 
In Figure~\ref{fig:j1_f} the outcome ($J_1$) of player~1 is depicted, for different strategies of the form $u_1(k) = K_1x(k)$, with the action of player~2 fixed at $u_2(k) = K_2^{OL} x(k)$.
It can be clearly appreciated, as discussed in Remark~\ref{rem:OL_vs_F} and Remark \ref{rem:comaprison} 
that the minimum of $J_1$ (indicated by the vertical red dashed line) is not, in general, attained by the OL-NE (indicated by the vertical black dotted line). A similar comparison is depicted in Figure~\ref{fig:j2_f}, by varying the strategy $u_2(k)= K_2x(k)$, with the action of player~1 fixed at $u_1(k) = K_1^{OL} x(k)$. This observation highlights the fact that the \emph{feedback synthesis of an OL-NE} does not yield a F-NE, since it does not in general enforce the inequalities~\eqref{eq:def_FNE} whenever the feasible set of alternative strategies includes control actions fed back from the current value of the state.
Nonetheless, the outcome of player~1 for different strategies $u_1(k)=K_1x(k)$ is shown in Figure~\ref{fig:j1_ol} when the strategy of player~2 is fixed at\footnote{Note that $(A_{cl}^{OL})^{k}x(0)$ corresponds to the trajectory $x(k)$ of the closed-loop system \eqref{eq:lti}, \eqref{eq:u_OL_NE_feedback_synth}, $i=1,2$.} $u_2(k) = K_2^{OL}(A_{cl}^{OL})^kx(0)$, where $A_{cl}^{OL}=\left(A + B_1K_1^{OL} + B_2K_2^{OL}\right)$. Then, the individual optimal control problem for the player 1 is solved by $u_1(k) = K_1^{OL}x(k)$, even if the admissible set includes feedback control laws, as can be clearly appreciated by the (intersecting) vertical red (indicating the minimum of $J_1$) and black dashed lines (indicating the outcome corresponding to $u_1(k) = K_1^{OL}x(k)$). The same observation can be made for player~2, as depicted in Figure~\ref{fig:j2_ol}. 
Note that in Figures \ref{fig:j1_f}-\ref{fig:j2_ol} the red stars represent the minimum of the cost functional.

\begin{figure}[ht]
\begin{minipage}{.5\textwidth}
\begin{center}
\includegraphics[width=0.8\columnwidth]{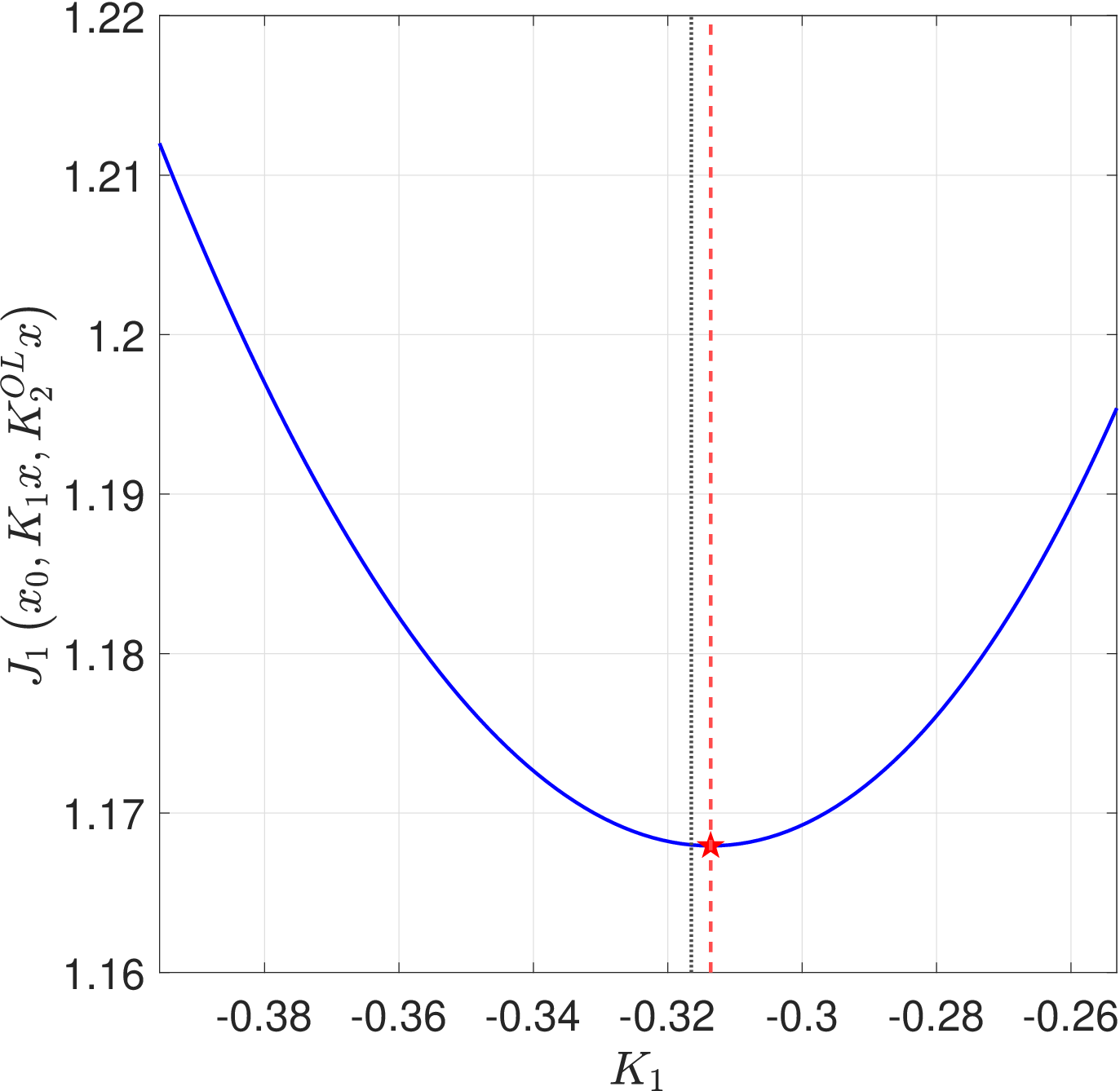}
\captionsetup{justification=centering,margin=0.1cm}
\caption{Graph of $J_1$ (blue line)
		corresponding to different strategies 
		$u_1(k) = K_1x(k)$ 
		with $u_2(k) = K_2^{OL} x(k)$. The vertical dotted line (black) corresponds to the gain $K_1^{OL}$ while the vertical dashed line (red) represents $\arg \min_{K_1} J_1(x_0,K_1x(k), K_2^{OL}x(k))$.}
  \label{fig:j1_f}
\end{center}\end{minipage}%
\begin{minipage}{.5\textwidth}
\begin{center}
\includegraphics[width=0.8\columnwidth]{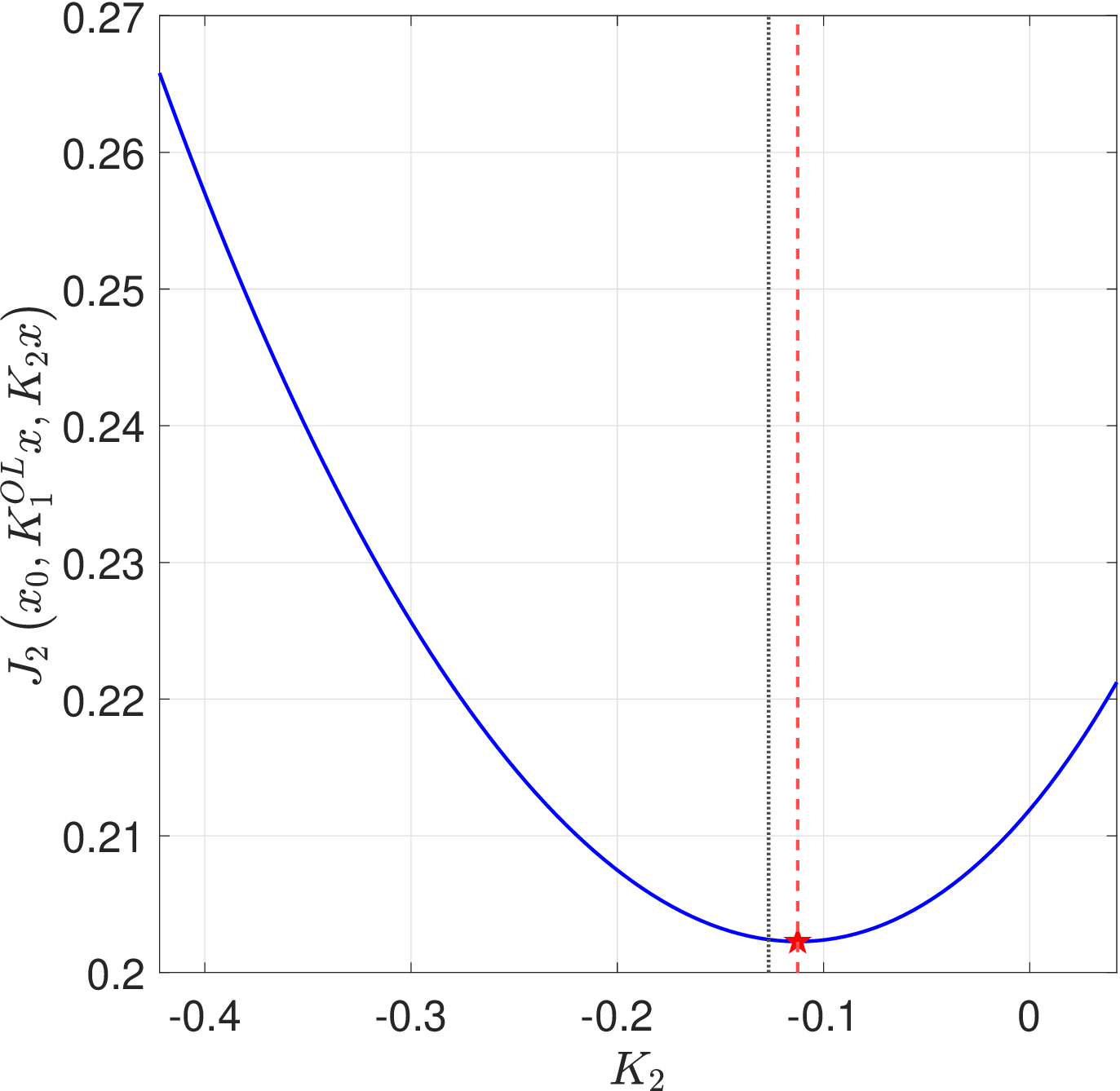}
\captionsetup{justification=centering,margin=0.1cm}
\caption{Graph of $J_2$ (blue line)
		corresponding to 
		different strategies 
		$u_2(k) = K_2x(k)$ 
		with $u_1(k) = K_1^{OL} x(k)$. The vertical dotted line (black) corresponds to the gain $K_2^{OL}$ while the vertical dashed line (red) represents $\arg \min_{K_2} J_2(x_0, K_1^{OL}x(k), K_2x(k))$.}
  \label{fig:j2_f}
\end{center}
\end{minipage}
\end{figure}

\begin{figure}[ht]
\begin{minipage}{.5\textwidth}
\begin{center}
\includegraphics[width=0.8\columnwidth]{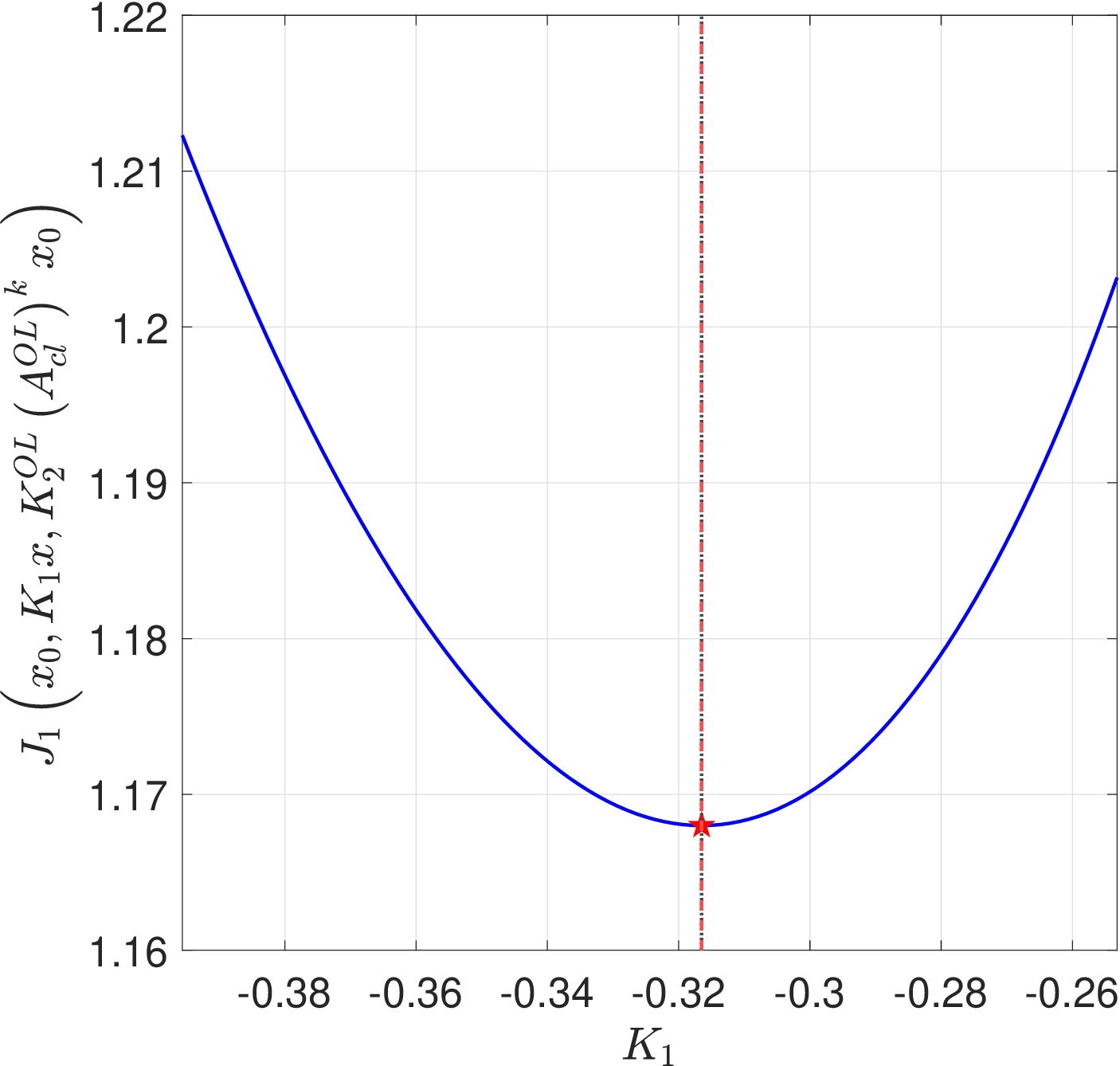}
\captionsetup{justification=centering,margin=0.1cm}
\caption{Graph of $J_1$ (blue line)
		corresponding to 
		different strategies 
		$u_1(k) = K_1x(k)$ 
		with $u_2(k) = K_2^{OL}\left(A_{cl}^{OL}\right)^{k}x(0)$. The vertical dotted line (black) corresponds to the gain $K_1^{OL}$ and coincides with the vertical dashed line (red) representing $\arg \min_{K_1} J_2(x_0,K_1x(k),\allowbreak K_2^{OL}(A_{cl}^{OL})^kx(0))$.}
  \label{fig:j1_ol}
\end{center}\end{minipage}%
\begin{minipage}{.5\textwidth}
\begin{center}
\includegraphics[width=0.8\columnwidth]{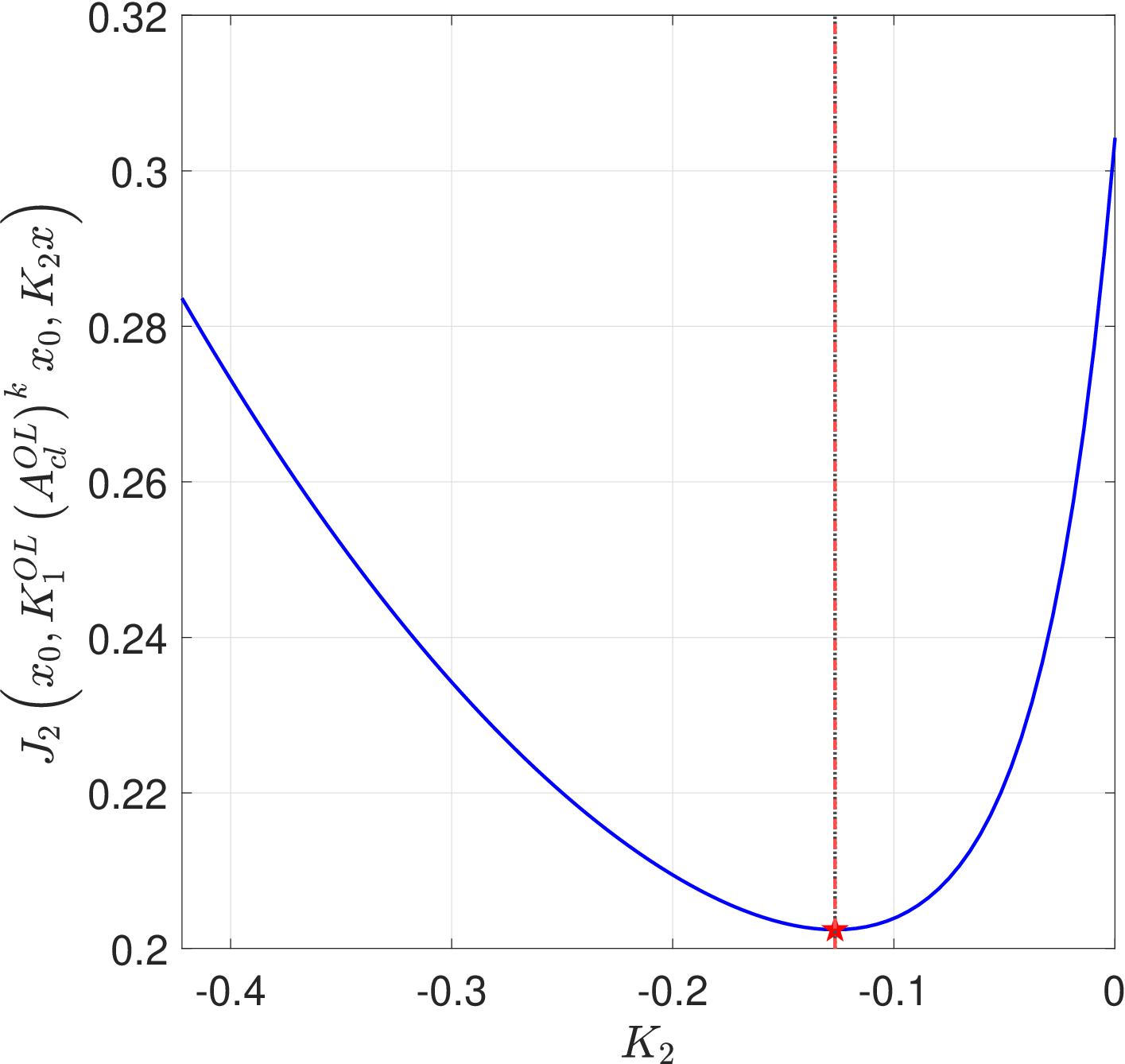}
\captionsetup{justification=centering,margin=0.1cm}
\caption{Graph of $J_2$ (blue line)
		corresponding to 
		different strategies 
		$u_2(k) = K_2x(k)$ 
		with $u_1(k) = K_1^{OL}\left(A_{cl}^{OL}\right)^{k}x(0)$. The vertical dotted line (black) corresponds to the gain $K_2^{OL}$ and coincides with the vertical dashed line (red) representing $\arg \min_{K_2} J_2(x_0, K_2^{OL}(A_{cl}^{OL})^kx(0), \allowbreak K_2x(k))$.}
  \label{fig:j2_ol}
\end{center}
\end{minipage}
\end{figure}

\addtolength{\textheight}{-0.5cm}
\section{Conclusions \hlb{and further extensions}}\label{sec:conc}
Feedback and Open-Loop NE have been studied and discussed for dynamic games defined by linear, discrete-time dynamics and quadratic cost functionals over an infinite horizon. Leveraging on DP and PMP, respectively, for the former we provide a comprehensive characterization in terms of the stabilizing solutions of coupled algebraic matrix equations, 
whereas for the latter 
we have derived 
an entirely new characterization in terms of certain quadratic matrix equations. This result has been enabled by an extension of PMP arguments to the context of infinite-horizon optimal control problems, with the underlying dynamics influenced by a known exogenous input. 
The analysis provides insights that are demonstrated and visualized on a numerical example. \hlb{As a further extension in the setting of OL-NE, it would be particularly relevant to relax the invertibility assumption imposed on the (open-loop) dynamic matrix, which appears crucial  in 
the arguments employed in this work. Furthermore,  considering practical motivations, 
it may be of interest to envision strategies that permit the construction of the equilibrium strategies discussed above on the basis of collected data only, thus establishing data-driven implementations of the corresponding results.} \color{black}

\bibliographystyle{siamplain}
\bibliography{biblio}

\begin{thebibliography}{10}

\bibitem{basar2008h}
{\sc T.~Ba{\c{s}}ar and P.~Bernhard}, {\em H-infinity optimal control and related minimax design problems: a dynamic game approach}, Springer Science \& Business Media, 2008.

\bibitem{basar1998}
{\sc T.~Ba{\c{s}}ar and G.~J. Olsder}, {\em Dynamic noncooperative game theory}, SIAM, 1998.

\bibitem{bertsekas2012dynamic}
{\sc D.~Bertsekas}, {\em Dynamic programming and optimal control: Volume II}, vol.~4, Athena scientific, 2012.

\bibitem{Cappello2021}
{\sc D.~Cappello, S.~Garcin, Z.~Mao, M.~Sassano, A.~Paranjape, and T.~Mylvaganam}, {\em A hybrid controller for multi-agent collision avoidance via a differential game formulation}, IEEE Transactions on Control Systems Technology, 29 (2021), pp.~1750--1757.

\bibitem{Cappello2022}
{\sc D.~Cappello and T.~Mylvaganam}, {\em Distributed differential games for control of multi-agent systems}, IEEE Transactions on Control of Network Systems, 9 (2022), pp.~635--646.

\bibitem{chen2019control}
{\sc F.~Chen, W.~Ren, et~al.}, {\em On the control of multi-agent systems: A survey}, Foundations and Trends{\textregistered} in Systems and Control, 6 (2019), pp.~339--499.

\bibitem{engwerda_infinityHor}
{\sc J.~Engwerda}, {\em The solution of the infinite horizon tracking problem for discrete time systems possessing an exogenous component}, Journal of Economic Dynamics and Control, 14 (1990), pp.~741--762.

\bibitem{Engwerda2005}
{\sc J.~Engwerda}, {\em LQ Dynamic Optimization and Differential Games}, Johm Wiley \& Sons, 2005.

\bibitem{Engwerda2007}
{\sc J.~Engwerda}, {\em Algorithms for computing nash equilibria in deterministic {LQ} games}, Computational Management Science, 4 (2007), pp.~113--140.

\bibitem{basar2019}
{\sc S.~R. Etesami and T.~Başar}, {\em Dynamic games in cyber-physical security: An overview}, Dynamic Games and Applications, 9 (2019), pp.~884--913.

\bibitem{freiling1999discrete}
{\sc G.~Freiling, G.~Jank, and H.~Abou-Kandil}, {\em Discrete-{T}ime {R}iccati {E}quations in {O}pen-{L}oop {N}ash and {S}tackelberg {G}ames}, European journal of control, 5 (1999), pp.~56--66.

\bibitem{golub}
{\sc G.~H. Golub and C.~F. Van~Loan}, {\em Matrix computations}, The John Hopkins University Press, 3rd~ed., 1996.

\bibitem{jaskiewicz2011stochastic}
{\sc A.~Ja{\'s}kiewicz and A.~S. Nowak}, {\em Stochastic games with unbounded payoffs: applications to robust control in economics}, Dynamic Games and Applications, 1 (2011), pp.~253--279.

\bibitem{jaskiewicz2022}
{\sc A.~Jaśkiewicz and A.~S. Nowak}, {\em A note on topological aspects in dynamic games of resource extraction and economic growth theory}, Games and Economic Behavior, 131 (2022), pp.~264--274.

\bibitem{jorgensen2010}
{\sc S.~Jørgensen, G.~Martín-Herrán, and G.~Zaccour}, {\em Dynamic games in the economics and management of pollution}, Environmental Modeling \& Assessment, 15 (2010), pp.~433--467.

\bibitem{kordonis2022}
{\sc I.~Kordonis, A.-R. Lagos, and G.~P. Papavassilopoulos}, {\em Dynamic games of social distancing during an epidemic: Analysis of asymmetric solutions}, Dynamic Games and Applications, 12 (2022), pp.~214--236.

\bibitem{kremer}
{\sc D.~Kremer}, {\em Non-Symmetric {R}iccati Theory and Noncooperative Games}, PhD. Thesis RWTH Aachen. Aachener Beitr¨age zur Mathematik, Band 30, ISBN: 3-86073-632-9., 2003.

\bibitem{kuvcera1972discrete}
{\sc V.~Ku{\v{c}}era}, {\em The discrete {R}iccati equation of optimal control}, Kybernetika, 8 (1972), pp.~430--447.

\bibitem{Li-Bur:19}
{\sc Y.~Li, G.~Carboni, F.~Gonzalez, D.~Campolo, and E.~Burdet}, {\em Differential game theory for versatile physical human--robot interaction}, Nature Machine Intelligence, 1 (2019), pp.~36--43.

\bibitem{limebeer1992mixed}
{\sc D.~Limebeer, B.~Anderson, and B.~Hendel}, {\em Mixed ${H}_2$/${H}_\infty$ filtering by the theory of {N}ash games}, in Robust Control, Springer, 1992, pp.~9--15.

\bibitem{huang2020}
{\sc H.~Linan and Z.~Quanyan}, {\em A dynamic games approach to proactive defense strategies against advanced persistent threats in cyber-physical systems}, Computers \& Security, 89 (2020), p.~101660.

\bibitem{long2011}
{\sc N.~V. Long}, {\em Dynamic games in the economics of natural resources: A survey}, Dynamic Games and Applications, 1 (2011), pp.~115--148.

\bibitem{Myl2017}
{\sc T.~Mylvaganam, M.~Sassano, and A.~Astolfi}, {\em A differential game approach to multi-agent collision avoidance}, IEEE Transactions on Automatic Control, 62 (2017), pp.~4229--4235.

\bibitem{nadini2020}
{\sc M.~Nadini, L.~Zino, A.~Rizzo, and M.~Porfiri}, {\em A multi-agent model to study epidemic spreading and vaccination strategies in an urban-like environment}, Applied Network Science, 5 (2020), pp.~68--97.

\bibitem{nagata2002}
{\sc T.~Nagata and H.~Sasaki}, {\em A multi-agent approach to power system restoration}, IEEE Transactions on Power Systems, 17 (2002), pp.~457--462.

\bibitem{Nor-Myl:ECC22}
{\sc B.~Nortmann and T.~Mylvaganam}, {\em Data-driven cost representation for optimal control and its relevance to a class of asymmetric linear quadratic dynamic games}, in 2022 European Control Conference, 2022, pp.~2185--2190.

\bibitem{sassano2022analysis}
{\sc M.~Sassano, T.~Mylvaganam, and A.~Astolfi}, {\em On the analysis of open-loop nash equilibria admitting a feedback synthesis in nonlinear differential games}, Automatica, 142 (2022), p.~110389.

\bibitem{pratap2017}
{\sc V.~P. Singh, N.~Kishor, and P.~Samuel}, {\em Distributed multi-agent system-based load frequency control for multi-area power system in smart grid}, IEEE Transactions on Industrial Electronics, 64 (2017), pp.~5151--5160.

\bibitem{starr1969further}
{\sc A.~W. Starr and Y.-C. Ho}, {\em Further properties of nonzero-sum differential games}, Journal of Optimization Theory and Applications, 3 (1969), pp.~207--219.

\bibitem{Starr:69}
{\sc A.~W. Starr and Y.-C. Ho}, {\em Nonzero-sum differential games}, Journal of optimization theory and applications, 3 (1969), pp.~184--206.

\bibitem{sweriduk1994differential}
{\sc G.~Sweriduk and A.~Calise}, {\em A differential game approach to the mixed ${H}_2$/${H}_\infty$ problem}, in Guidance, Navigation, and Control Conference, 1994, pp.~1072--1082.

\bibitem{Wrz20}
{\sc M.~Wrzos-Kaminska, T.~Mylvaganam, K.~Y. Pettersen, and J.~Tommy~Gravdahl}, {\em Collision avoidance using mixed {$H_2$/$H_\infty$} control for an articulated intervention-{AUV}}, in 2020 European Control Conference, 2020, pp.~881--888.

\end{thebibliography}
\end{document}